\documentclass[12pt]{article}
\usepackage{amsthm,amsfonts,amssymb,amscd}

\begin{document}

\begin{center}
\large \bf Birational geometry of singular \\Fano hypersurfaces
of index two
\end{center}\vspace{0.5cm}

\centerline{A.V.Pukhlikov}\vspace{0.5cm}

\parshape=1
3cm 10cm \noindent {\small \quad\quad\quad \quad\quad\quad\quad
\quad\quad\quad {\bf }\newline For a Zariski general
(regular) hypersurface $V$ of degree $M$ in the $(M+1)$-dimensional
projective space, where $M\geqslant 16$, with at most quadratic
singularities of rank $\geqslant 13$, we give a complete description
of the structures of rationally connected (or Fano-Mori) fibre space:
every such structure over a positive-dimensional base is a pencil
of hyperplane sections. This implies, in particular, that $V$ is non-rational and its groups of birational and biregular automorphisms coincide: $\mathop{\rm Bir} V = \mathop{\rm Aut} V$. The set of non-regular hypersurfaces has codimension at least $\frac12(M-11)(M-10)-10$ in the natural parameter space.

Bibliography: 25 titles.} \vspace{1cm}

\vspace{1cm}

\noindent 14E05, 14E07\vspace{1cm}

\section*{Introduction}

{\bf 0.1. Statement of the main result.} The symbol ${\mathbb
P}$ stands for the complex projective space ${\mathbb
P}^{M+1}$, $M\geqslant 4$. Hypersurfaces of degree $M$ in
${\mathbb P}$ are parametrized by the points of the projective space
$$
{\cal F}={\mathbb P}(H^0({\mathbb P},{\cal O}_{{\mathbb P}}(M))).
$$
Let $V\subset{\mathbb P}$ be a hypersurface of degree $M$. If it is irreducible, reduced, factorial and has at most terminal singularities, then $V$ is a Fano variety of index two:
$$
\mathop{\rm Pic} V= \mathop{\rm Cl} V= {\mathbb Z} H,\quad K_V=-2H,
$$
where $H$ is the class of hyperplane section. If $P\subset {\mathbb
P}$ is a linear subspace of codimension 2, then restricting to the hypersurface $V$ the linear projection $\alpha_{P}\colon
{\mathbb P}\dashrightarrow {\mathbb P}^1$ from the subspace $P$, we define on $V$ a structure of a Fano-Mori fibre space $\pi_P\colon
V\dashrightarrow {\mathbb P}^1$.\vspace{0.1cm}

Let $\lambda\colon Y\to S$ be a rationally connected fibre space, that is, a surjective morphism of projective varieties,
$\mathop{\rm dim} S\geqslant 1$, where the fibre of general position
$\lambda^{-1}(s)$, $s\in S$, and the base $S$ are rationally connected.\vspace{0.1cm}

Here is the main result of the present paper.\vspace{0.1cm}

{\bf Theorem 1.} {\it For $M\geqslant 16$ there is a Zariski open subset ${\cal U}\subset {\cal F}$, such that:\vspace{0.1cm}

{\rm (i)} every hypersurface $V\in {\cal U}$ is irreducible, reduced, factorial and has at most terminal singularities;\vspace{0.1cm}

{\rm (ii)} the inequality
$$
\mathop{\rm codim}(({\cal F}\setminus {\cal U})\subset
{\cal F})\geqslant \frac12(M-11)(M-10)-10
$$
holds,\vspace{0.1cm}

{\rm (iii)} for every $V\in{\cal U}$ and every birational map $\chi\colon V\to Y$ onto the total space of the rationally connected fibre space $\lambda\colon Y\to S$ over a positive-dimensional base $S$ we have $S={\mathbb P}^1$ and for some isomorphism $\beta\colon {\mathbb P}^1\to S$ and some subspace $P\subset {\mathbb P}$ of codimension 2 the equality
$$
\lambda\circ\chi=\beta\circ\pi_P,
$$
holds, that is, the following diagram commutes:}
$$
\begin{array}{ccccc}
& V & \stackrel{\chi}{\dashrightarrow} & Y &\\
\pi_P\!\!\!\! & \downarrow & & \downarrow &\!\!\!\! \lambda\\
& {\mathbb P}^1 & \stackrel{\beta}{\to} & S. &\\
\end{array}
$$

For a Zariski general smooth hypersurface $V\subset{\mathbb
P}$ the claim (iii) of Theorem 1 was shown for $M\geqslant 16$ in
\cite{Pukh16a}. Hypersurfaces with at least one singular point form a divisor in the space ${\cal F}$. Thus in the present paper we essentially improve the main result of \cite{Pukh16a}: we extend it to hypersurfaces with bounded singularities and give an effective estimate for the codimension of the complement to the set ${\cal U}$ of ``correct'' hypersurfaces
(which grows as $\frac12 M^2$ when the dimension $M$ grows).\vspace{0.1cm}

Theorem 1 immediately implies the standard set of facts about birational geometry of the variety $V\in{\cal U}$.\vspace{0.1cm}

{\bf Corollary 1.} {\it For every hypersurface $V\in{\cal U}$
the following claims are true.\vspace{0.1cm}

{\rm (i)} On the variety $V$ there are no structures of a rationally connected fibre space (and therefore, of a Fano-Mori fibre space) over a base of dimension $\geqslant 2$. In particular, on $V$ there are no structures of a conic bundle and del Pezzo fibrations, and the variety $V$ itself is non-rational.\vspace{0.1cm}

{\rm (ii)} Assume that there is a birational map
$\chi\colon V\dashrightarrow Y$, where $Y$ is a Fano variety of index $r\geqslant 2$ with factorial terminal singularities, such that $\mathop{\rm Pic}Y={\mathbb Z}H_Y$, where $K_Y=-rH_Y$, and the linear system $|H_Y|$ is non-empty and free. Then $r=2$ and the map $\chi$ is abiregular isomorphism.\vspace{0.1cm}

{\rm (iii)} The groups of biregular and birational automorphisms of the variety $V$ are equal: $\mathop{\rm Bir}V=\mathop{\rm Aut}V$.
}\vspace{0.3cm}


{\bf 0.2. The regularity conditions.} Now we give an explicit description of the open set ${\cal F}_{\rm reg}\subset{\cal F}$, consisting of hypersurfaces, satisfying the regularity conditions, stated below. We will show that for ${\cal U}={\cal F}_{\rm reg}$ all statements of Theorem 1 are true.\vspace{0.1cm}

Let $o\in{\mathbb P}$ be an arbitrary point, $(z-1,\dots,
z_{M+1})=(z_*)$ a system of affine coordinates with the origin at the point
$o$ and $V\ni o$ a hypersurface of degree $M$. It is defined by an equation $f=0$, where
$$
f=q_1+q_2+\cdots +q_M
$$
is a non-homogeneous polynomial in the variables $z_*$, decomposed into homogeneous components $q_i$ of degree $i\geqslant 1$. The regularity conditions depend on whether the point $o\in V$ is singular or non-singular, that is, whether $q_1\equiv 0$ èëè $q_1\not\equiv 0$.\vspace{0.1cm}

First, we state the regularity conditions for a {\bf non-singular point}.\vspace{0.1cm}

(R1.1) For every linear subspace $\Pi\subset{\mathbb
C}^{M+1}$ of the standard coordinate space with the coordinates
$z_*$, of codimension $\mathop{\rm codim}\Pi=c\in\{0,1,2,3\}$ and such that $q_1|_{\Pi}\not\equiv 0$ (that is, $\Pi\not\subset T_oV$), the sequence
$$
q_1|_{\Pi}, \,\,q_2|_{\Pi},\,\,\dots, \,\,q_{M-c}|_{\Pi}
$$
is regular in the local ring ${\cal O}_{o,\Pi}$, that is, the system of equations
$$
q_1|_{\Pi}, \,=\,q_2|_{\Pi},\,=\,\dots, \,=\,q_{M-c}|_{\Pi}=0
$$
defines a finite set of lines through the point $o$ (and in the projective space ${\mathbb P}(\Pi)$ a finite set of points).\vspace{0.1cm}

(R1.2) The rank of the quadratic form
$$
q_2|_{\{q_1=0\}}=q_2|_{T_oV}
$$
is at least 11.\vspace{0.1cm}

(R1.3) For every linear subspace $\Lambda\subset T_oV$ of codimension 2 the system of equations
$$
q_2|_{\Lambda}=q_3|_{\Lambda}=0
$$
defines an irreducible reduced closed set of codimension
2 in $\Lambda$.\vspace{0.1cm}

Now let us consider the regularity conditions for a {\bf singular point},
that is to say, we assume that $q_1\equiv 0$.\vspace{0.1cm}

(R2.1) For every linear subspace $\Pi\subset{\mathbb
C}^{M+1}$ of codimension $c\in\{1,2,3\}$ the sequence
$$
q_2|_{\Pi}, \,\,q_3|_{\Pi},\,\,\dots, \,\,q_{M+1-c}|_{\Pi}
$$
is regular in ${\cal O}_{o,\Pi}$, that is, the system of equations
$$
q_2|_{\Pi}, \,=\,q_3|_{\Pi},\,=\,\dots, \,=\,q_{M+1-c}|_{\Pi}=0
$$
defines a finite set of lines through the point
$o$.\vspace{0.1cm}

(R2.2) The rank of the quadratic form $q_2$ is at least 13.\vspace{0.1cm}

(R2.3) For every linear subspace $\Pi\subset{\mathbb
C}^{M+1}$ of codimension 2 the closed set
$$
\{q_2|_{{\mathbb P}(\Pi)}=q_3|_{{\mathbb P}(\Pi)}=0\}
$$
in the projective space ${\mathbb P}(\Pi)\cong{\mathbb
P}^{M-2}$ is irreducible, reduced factorial complete intersection of type $2\cdot 3$.\vspace{0.1cm}

{\bf Definition 1.} We say that the hypersurface $V\in {\cal F}$ is
{\it regular}, if at every non-singular point $o\in V$ it is regular in the sense of the conditions (R1.1-3), and in every singular point $o\in
V$ it is regular in the sense of the conditions (R2.1-3).\vspace{0.1cm}

The set of regular hypersurfaces is denoted by the symbol ${\cal
F}_{\rm reg}$. Obviously, ${\cal F}_{\rm reg}\subset{\cal F}$ is a Zariski open subset.\vspace{0.1cm}

The condition (R2.2) implies that the codimension of the singular set $\mathop{\rm Sing} V$ of a regular hypersurface $V$ is at least 10, so that $V$ is irreducible, reduced and by the well known theorem of Grothendieck \cite{CL}, factorial. The same condition (R2.2) implies that the singularities of a regular hypersurface $V$ are terminal (see \cite{EP} and also \cite{Pukh15a}; in the latter paper at the end of Subsection 2.1 it is  explained that quadratic singularities, the rank of which is bounded from below, are stable with respect to blow ups, which, in its turn, makes it very easy to see the property of being terminal). Therefore, for a regular hypersurface $V$ the claim (i) of Theorem 1 is true. By what was said, Theorem 1 is implied by the following two facts.\vspace{0.1cm}

{\bf Theorem 2.} {\it The open set ${\cal U}={\cal
F}_{\rm reg}$ satisfies the claim (ii) of Theorem 1.}\vspace{0.1cm}

{\bf Theorem 3.} {\it Every regular hypersurface
$V\in{\cal F}_{\rm reg}$ satisfies the claim (iii) of Theorem
1.}\vspace{0.3cm}


{\bf 0.3. The method of maximal singularities.} For an arbitrary subspace $P\subset{\mathbb P}$ of codimension 2 by the symbol
$|H-P|$ we denote the pencil of divisors cut out on $V$ by the pencil of hyperplanes containing $P$. In the notations of the part (iii) of Theorem 1 let $\Sigma_Y$ be the $\lambda$-pull back of some very ample linear system on the base $S$, and $\Sigma$ its strict transform on $V$ with respect to $\chi$. The linear system $\Sigma$ is mobile (that is, has no fixed components) and we may assume that for some $n\geqslant 1$
$$
\Sigma\subset |2nH|
$$
(replacing, if necessary, the very ample system on the base $S$ by its symmetric square). This whole set of geometric objects: the hypersurface $V\in{\cal F}_{\rm reg}$, the rationally connected fibre space $\lambda\colon Y\to S$, the birational map $\chi$, the linear systems $\Sigma_Y$ on $Y$ and $\Sigma$ on $V$, and therefore, the number $n\geqslant 1$, is assumed to be fixed.\vspace{0.1cm}

It is well known (see, for instance, \cite[Chapter 2, Section
1]{Pukh13a}, and also Subsection 1.1 of the present paper), that the mobile linear system $\Sigma$ has a {\it maximal singularity}: for some exceptional divisor $E^*$ over $V$ the {\it Noether-Fano inequality} holds:
\begin{equation}\label{23.10.2017.1}
\mathop{\rm ord}\nolimits_{E^*}\Sigma>n\cdot a(E^*),
\end{equation}
where $a(E^*)$ is the discrepancy of $E^*$ with respect to the model $V$. In a different way this can be expressed as follows: the pair $(V,\frac{1}{n}D)$ is not canonical for a general divisor $D\in\Sigma$ or, even simpler, the pair $(V,\frac{1}{n}\Sigma)$ is not canonical.\vspace{0.1cm}

There are mobile linear systems with a maximal singularity on $V$. For instance, let $E_P$ be the exceptional divisor of the blow up of the subvariety $V\cap P$ of codimension 2 on $V$, where $P\subset{\mathbb P}$ is a linear subspace of codimension 2. Obviously, $a(E_P)=1$, so that the  ``double pencil'' $|H-P|$, that is, the linear system $|2H-2P|\subset|2H|$,  has $E_P$ as a maximal singularity, since
$$
\mathop{\rm ord}\nolimits_{E_P}|2H-2P|=2.
$$

Theorem 3 essentially means that {\it any} linear system
$\Sigma$ with a maximal singularity is composed of a pencil
$|H-P|$. The proof of Theorem 3 consists of two main steps.\vspace{0.1cm}

{\bf Theorem 4.} {\it Assume that for a certain linear subspace $P\subset{\mathbb P}$ of codimension 2 the inequality
\begin{equation}\label{21.10.2017.1}
\mathop{\rm mult}\nolimits_{P\cap V}\Sigma>n
\end{equation}
holds. Then $\Sigma$ is composed of the pencil $|H-P|$, that is, every divisor $D\in\Sigma$ is a sum of $2n$ hyperplane sections from that pencil.}\vspace{0.1cm}

{\bf Theorem 5.} {\it For a linear system $\Sigma$ with a maximal singularity there is a linear subspace $P\subset{\mathbb P}$ of codimension 2 satisfying the inequality (\ref{21.10.2017.1}).}\vspace{0.1cm}

Theorem 3 obviously follows from Theorems 4 and 5. Proof of Theorem 5 is most difficult. Proof of Theorem 2 is not hard.\vspace{0.3cm}


{\bf 0.4. The structure of the paper.} In \S 1 we show Theorem 4. The following fact is crucial in the proof of Theorem 4: the global log canonical threshold of {\it every} hyperplane section of the hypersurface
$V$ is equal to 1. The equality $\mathop{\rm lct}
(F)=1$ for Fano hypersurfaces $F\subset{\mathbb P}^M$ of degree
$M$, satisfying certain restrictions for the singularities and the regularity conditions at non-singular and singular points, has been recently proven in \cite{Pukh15a}, so that in this paper we just check that every hyperplane section of the hypersurface $V$ satisfies the requirements of \cite{Pukh15a}.\vspace{0.1cm}

In \S 2 we prove Theorem 2. For each of the regularity conditions we estimate the codimension of the set of hypersurfaces which do not satisfy that condition at at least one point. After that, by means of the technique of hypertangent divisors we prove certain estimates, bounding the multiplicities of irreducible subvarieties of the hypersurface $V$ at singular points $o\in\mathop{\rm Sing} V$ from below. Those estimates will be needed later.\vspace{0.1cm}

In \S 3 we start the proof of Theorem 5. Following the traditional scheme of arguments of the method of maximal singularities, we assume that there is no linear subspace $P\subset{\mathbb P}$ of codimension 2, satisfying the inequality (\ref{21.10.2017.1}). We have to show that under this assumption the linear system $\Sigma$ has no maximal singularities at all: this contradiction proves Theorem 5. In Subsection 3.1 we prove that the centre
$B^*$ of the maximal singularity $E^*$ is contained in the singular locus
$\mathop{\rm Sing} V$ of the hypersurface $V$. In order to do this, it is sufficient to check that if $B\not\subset\mathop{\rm Sing} V$, then the maximal singularity $E^*$ is excluded by the arguments of
\cite{Pukh16a}.\vspace{0.1cm}

In Subsection 3.2 for a point $o\in B$ of general position (which by what we have already proven is a quadratic singularity of the hypersurface
$V$) we prove that there is, generally speaking, another singularity $E$ of the linear system $\Sigma$, satisfying a {\it Noether-Fano type inequality}, which is weaker than (\ref{23.10.2017.1}), but still strong enough for our purposes. We do it by means of the inversion of adjunction \cite{Kol93} similar to the arguments of \cite[Subsection 4.2]{Pukh16a}. Finally, in Subsection 3.3 we recall the improved version of the technique of counting multiplicities for a complete intersection singularity \cite{Pukh2017a}.\vspace{0.1cm}

In \S 4 we prove certain technical statements about the secant variety of a subvariety of small codimension on a quadratic hypersurface of sufficiently high rank; those technical facts are used in \S 5 for exclusion of the singularity $E$.\vspace{0.1cm}

\S 5 is the central part of the proof of Theorem 5. Depending on the type of the singularity $E$ (the types are defined in Subsection 3.2), it is excluded by different methods. In accordance with the traditional scheme of the method of maximal singularities (see \cite[Chapter 2]{Pukh13a}), we consider the {\it self-intersection} $Z=(D_1\circ D_2)$ of the mobile linear system $\Sigma$ and prove that the existence of the singularity $E$ imposes so strong restrictions on the singularities of the self-intersection $Z$, which can not be satisfied for an effective cycle of codimension 2 on $V$. Thus we prove that the mobile linear system $\Sigma$ can not have the singularity $E$, and therefore can not have the maximal singularity $E^*$, either. This contradiction completes the proof of Theorem 5 (and the main Theorem 1).\vspace{0.3cm}


{\bf 0.5. Historical remarks and acknowledgements.} The few attempts to study birational geometry of higher-dimensional Fano varieties of index higher than 1 are listed in the introduction to \cite{Pukh16a}, see also the introduction to \cite{Pukh15a}. Here we note that, starting from the paper \cite{EP}, the results about birational rigidity of particular classes of Fano varieties become effective in the sense that an explicit effective estimate for the codimension of the subset of non-rigid varieties in the natural parameter space of the given family is produced. Those results (see \cite{EvP2017,Johnstone}) are very important because they open the way for the study of the problem of birational rigidity for Fano fibre spaces over a base of high (ideally --- arbitrary) dimension, the fibres of which belong to a given family of Fano varieties. The first breakthrough in that direction is the paper \cite{Pukh15a}. In the present paper the result of \cite{Pukh16a} is extended to singular Fano hypersurfaces of index 2 and becomes effective in the sense described above: we give an explicit effective estimate of the set of hypersurfaces, the birational geometry of which does not satisfy the property (iii) of Theorem 1.\vspace{0.1cm}

Recently quite a few papers were produced, proving the stable non-rationality of various classes of Fano varieties and Fano-Mori fibre space, see, for instance, [11-19] (the list is by no means complete). The importance of those results, obtained by completely different methods (compared to the method of maximal singularities), can not be overestimated. Note, however, that the stable non-rationality is shown for a {\it very general} variety in the family. The method of maximal singularities gives birational rigidity (or an explicit exhaustive description of birational geometry like what is done in this paper) for a {\it Zariski general} variety, together with an effective estimate for the codimension of the complement in the parameter space.\vspace{0.1cm}

Note also the recent paper on the birational rigidity of singular Fano three-folds \cite{ChelPark2017}, the recently published paper \cite{ChelGrin2017} and the papers
\cite{ProkhShr2014,ProkhShr2016} on the groups of birational automorphisms.\vspace{0.1cm}

The author thanks The Leverhulme Trust for the financial support of the present project (Research Project Grant RPG-2016-279).\vspace{0.1cm}

The author is also grateful to the members of the Divisions of Algebraic Geometry and Algebra of Steklov Mathematical Institute for the interest to this work and also to the colleagues -- algebraic geometers at the University of Liverpool for the general support.


\section{The pencils of hyperplane section \\ and the regularity conditions}

In the section we prove Theorem 4. As the first step of the proof, we consider the new model of the hypersurface $V$, which is obtained by blowing up the subvariety $V\cap P$ (Subsection 1.1). After that we get the alternative: either the claim of Theorem 4 is satisfied, or the strict transform of the linear system $\Sigma$ on the new model again has a maximal singularity (Subsection 1.2). Finally, in Subsection 1.3 we show that the results of \cite{Pukh15a} imply that the second case does not realize, because the global log canonical threshold of every hyperplane section of the hypersurface $V$ is equal to 1. This completes the proof of Theorem 4.\vspace{0.3cm}

{\bf 1.1. The structure of a Fano fibre space.} Let us prove Theorem 4. Set
$B=P\cap V$. Let $\varphi\colon V^+\to V$ be the blow up of the subvariety $B$. Denote by the symbol $E_B$ the exceptional divisor of this blow up. The variety $V^+$ can be seen as the strict transform of the hypersurface $V$ with respect to the blow up
$\varphi_{{\mathbb P}}\colon {\mathbb P}^+\to{\mathbb P}$
of the linear subspace $P$, so that $E_B=V^+\cap E_P$, where
$E_P=\varphi^{-1}(P)$. The linear projection ${\mathbb
P}\dashrightarrow {\mathbb P}^1$ from the subspace $P$ extends to a ${\mathbb P}^M$-bundle $\pi_{{\mathbb
P}}\colon {\mathbb P}^+\to{\mathbb P}^1$. Set
$\pi=\pi_{{\mathbb P}}|_{V^+}\colon V^+\to{\mathbb
P}^1$.\vspace{0.1cm}

{\bf Proposition 1.1.} {\it {\rm (i)} The variety $V^+$ and every fibre $F_t=\pi^{-1}(t)$, $t\in{\mathbb P}^1$, are factorial and have at most terminal singularities. Every fibre $F_t$,
$t\in{\mathbb P}^1$, is a Fano variety.\vspace{0.1cm}

{\rm (ii)} The equalities
$$
\mathop{\rm Pic} V^+={\mathbb Z} H\oplus{\mathbb Z} E_B={\mathbb
Z} K^+\oplus {\mathbb Z} F
$$
hold, where $H=\varphi^* H$ for simplicity of notations, $K^+=K_{V^+}$ is the canonical class of the variety $V^+$, $F$ is the class of the fibre of the projection
$\pi$ and}
$$
K^+=-2H+E_B, \quad F=H-E_B.
$$

{\bf Proof.} The fibres of the projection $\pi_{{\mathbb P}}$ are isomorphic to hyperplanes (containing the subspace $P$) in ${\mathbb P}$, so that the fibres $F_t$ are isomorphic to the corresponding hyperplane sections of the hypersurface $V$, that is, to hypersurfaces of degree $M$ in
${\mathbb P}^M$. The conditions (R1.2) and (R2.2) imply that every hypersurface $F_t\subset{\mathbb P}^M$, $t\in{\mathbb P}^1$, gas at most quadratic singularities of rank at least 11. Therefore, the variety $V^+$ also has at most quadratic singularities of rank $\geqslant 11$. The claim (i) follows from here. The claim (ii) is checked by obvious computations. Q.E.D. for the proposition.\vspace{0.1cm}

Now let us consider the strict transform $\Sigma^+$ of the linear system
$\Sigma$ on $V^+$. This is a mobile linear system, and for some $m\in{\mathbb Z}_+$ and $l\in{\mathbb Z}$ we have the inclusion
$$
\Sigma^+\subset |-mK^++lF|.
$$
The formulas of part (ii) of Proposition 1.1 imply that
$m=2n-\mathop{\rm mmm}\nolimits_B\Sigma$ and $l=2(\mathop{\rm
mmm}\nolimits_B\Sigma - n)\geqslant 2$. If $m=0$, then the linear system $\Sigma^+$ is composed of the pencil $|F|$, so that the system
$\Sigma$ is composed of the pencil $|H-P|$, as Theorem 4 claims. Therefore let us assume that $m\geqslant 1$, and show that this assumption leads to a contradiction.\vspace{0.3cm}


{\bf 1.2. Maximal singularities of the system $\Sigma^+$.} The following claim is true.\vspace{0.1cm}

{\bf Proposition 1.2.} {\it The linear system $\Sigma^+$ has a maximal singularity: for some exceptional divisor $E^+$ over $V^+$ the Noether-Fano inequality
$$
\mathop{\rm ord}\nolimits_{E^+}\Sigma^+ > m\cdot a(E^+,V^+)
$$
holds, that is, for a general divisor $D^+\in \Sigma^+$ the pair
$(V^+,\frac{1}{m}D^+)$ is not canonical.}\vspace{0.1cm}

{\bf Proof.} This is a particular case of a general well known fact, see \cite[Chapter 2, Section 1]{Pukh13a}. For the convenience of the reader we give a sketch of a proof. Let ${\widetilde
V}\to V^+$ be the resolution of singularities of the birational map
$\chi\circ\varphi\colon V^+\dashrightarrow Y$, and $\widetilde
\Sigma$ the strict transform of the linear system $\Sigma^+$ on
${\widetilde V}$. Furthermore, let ${\cal E}$ be the set of all prime exceptional divisors of the resolution ${\widetilde V}\to
V^+$. Recall that ${\widetilde \Sigma}$ is the pull back of the free linear system $\Sigma_Y$ on ${\widetilde V}$. Since divisors of the system $\Sigma_Y$ by assumption are pulled back from the base $S$, and the general fibre of the projection $\lambda\colon Y\to S$ is rationally connected, for a general divisor ${\widetilde D}\in {\widetilde \Sigma}$ the class
$$
{\widetilde D}+m{\widetilde K}
$$
is not pseudoeffective (${\widetilde K}=K_{{\widetilde V}}$ for the brevity of writing). However,
$$
{\widetilde D}\sim -mK^++lF-\sum_{E\in {\cal E}}(\mathop{\rm
ord}\nolimits_E \Sigma^+)E
$$
and
$$
{\widetilde K}\sim K^++\sum_{E\in {\cal E}}a(E,V^+)E
$$
(for simplicity of notations the pull back of a divisor is denoted by the same symbol as the divisor itself), so that
$$
{\widetilde D}+m{\widetilde K}\sim lF+\sum_{E\in {\cal
E}}(ma(E,V^+)-\mathop{\rm ord}\nolimits_E \Sigma^+)E.
$$
We conclude that in the right hand side for at least one $E$ the corresponding coefficient is negative. Q.E.D. for the proposition.\vspace{0.1cm}

{\bf Remark 1.1.} In a similar way one proves that the original linear system $\Sigma$ has a maximal singularity, see Subsection 0.3.\vspace{0.1cm}

Let $R\subset V^+$ be the centre of the maximal singularity $E^+$ on
$V^+$, so that $\mathop{\rm codim}(R\subset V^+)\geqslant 2$. There are two options:\vspace{0.1cm}

\begin{itemize}
\item $R$ covers the base ${\mathbb P}^1$: $\pi(R)={\mathbb P}^1$,

\item $\pi(R)$ is a point on ${\mathbb P}^1$.

\end{itemize}
Assume that the first option takes place. Restricting the linear system $\Sigma^+$ onto the fibre $F_t$ of general position, we obtain a mobile linear system $\Sigma_t\subset |mH_t|$, where
$H_t$ is the class of a hyperplane section of $F_t\subset{\mathbb P}^M$,
and moreover, the pair $(F_t, \frac{1}{m}\Sigma_t)$ is not canonical (that is, $\Sigma_t$ has a maximal singularity). In \cite{EP} it was shown that this is impossible (under the weaker assumptions about the singularities of the hypersurface $F_t\subset{\mathbb P}^M$ and for weaker, than in the present paper, regularity conditions at every point). Therefore we may assume that the second option takes place: $R\subset F_t$ for some $t\in{\mathbb P}^1$. Somewhat abusing the notations, we write $F$ instead of
$F_t$. Since the linear system $\Sigma^+$ is mobile, it can be restricted onto $F$ and by inversion of adjunction \cite[Section
17.4]{Kol93} obtain an effective divisor $D_F\sim mH_F$, such that the pair $(F,\frac{1}{m}D_F)$ is not log canonical. However, this contradicts to the following fact.\vspace{0.1cm}

{\bf Proposition 1.3.} {\it For every divisor $\Delta\in |mH_F|$ the pair $(F,\frac{1}{m}\Delta)$ is log canonical.}\vspace{0.1cm}

{\bf Proof} is given below in Subsection 1.3.\vspace{0.1cm}

The contradiction obtained above shows that the case $m\geqslant 1$ is impossible. Proof of Theorem 4 is complete. Q.E.D.\vspace{0.3cm}


{\bf 1.3. The global log canonical threshold of a fiber.} The claim of Proposition 1.3 is shown in \cite[Theorem 1.4]{Pukh15a} under the assumption that the hypersurface $F\subset{\mathbb P}^M$ satisfies certain regularity conditions at every point $o\in
F$ (the conditions of the same type that the conditions (R$\alpha.\beta$) of the present paper). Therefore, in order to prove Proposition 1.3, it is sufficient to compare the conditions used in \cite{Pukh15a} with the conditions in Subsection 0.2 of the present paper and make sure that the latter are not weaker. In order to make the reading more convenient, we reproduce the regularity conditions from \cite{Pukh15a} below. To avoid any misunderstanding, the condition which in \cite{Pukh15a} has number (R$\alpha.\beta$) (for instance, (R2.1)), will be denoted by (${\rm R}^*$ $\alpha.\beta$).\vspace{0.1cm}

So let $F\subset{\mathbb P}^M$ be a hypersurface of degree
$M$, $o\in F$ an arbitrary point, $(u_1,\dots, u_M)$ a system of affine coordinates with the origin at the point $o$ and
$$
w=q_1^*+q_2^*+q_3^*+\cdots +q_M^*
$$
the affine equation of the hypersurface $F$ with respect to that system of coordinates, decomposed into homogeneous components. Here is the list of conditions, which should be satisfied for the hypersurface $F$ in
\cite{Pukh15a}. Let us first consider a non-singular point $o\in
F$.\vspace{0.1cm}

(${\rm R}^*$ 1.1) The sequence $q_1^*,\dots, q_{M-1}^*$ is regular in ${\cal O}_{o,{\mathbb P}^M}$.\vspace{0.1cm}

This condition is satisfied because it is a particular case of the condition
(R1.1) (for $c=1$).\vspace{0.1cm}

(${\rm R}^*$ 1.2) The quadratic form $q^*_2|_{\{q^*_1=0\}}$ is of rank $\geqslant 6$ and the linear span of every irreducible component of the closed set
\begin{equation}\label{25.10.2017.1}
\{q_1^*=q_2^*=q_3^*=0\}
\end{equation}
is the hyperplane $\{q^*_1=0\}$.\vspace{0.1cm}

The first part of this condition follows from (R1.2) (the rank of the quadratic form $q^*_2|_{\{q^*_1=0\}}$ turns out to be at least 9), and the second part follows from (R1.3), since by the condition (R1.3) the closed set (\ref{25.10.2017.1}) is irreducible, reduced and of codimension 2 in the hyperplane $\{q_1^*=0\}$, that is, forms an irreducible and reduced complete intersection of type $2\cdot 3$.\vspace{0.1cm}

(${\rm R}^*$ 1.3) For every hyperplane $P\subset{\mathbb
P}^M$, $P\ni o$, $P\neq T_oF$, the algebraic cycle of the scheme-theoretic intersection
$$
(P\circ \overline{\{q^*_1=0\}}\circ \overline{\{q^*_2=0\}}\circ F)
$$
is irreducible and reduced.\vspace{0.1cm}

This condition holds in our case: a section of $F$ by two hyperplanes is a section of $V$ by three hyperplanes, and therefore it has at most quadratic singularities of rank $\geqslant 7$ and for that reason, $(P\circ \overline{\{q^*_1=0\}}\circ F)$ is a factorial hypersurface in the projective space $P\cap \overline{\{q^*_1=0\}}$. The restriction of the quadratic form $q^*_2$ (that is, the restriction of $q_2$) onto this projective space is of rank $\geqslant 7$. Therefore, the condition (${\rm R}^*$ 1.3) holds.\vspace{0.1cm}

Now let us consider a singular point $o\in F$. Çäåñü $q_1^*\equiv 0$, so that the equation $w$ starts with $q^*_2$.\vspace{0.1cm}

(${\rm R}^*$ 2.1) For every linear subspace
$\Pi_{\sharp}\subset {\mathbb C}^M$ of codimension $c\in\{0,1,2\}$
the sequence
$$
q^*_2|_{\Pi_{\sharp}}, \dots, q^*_{M-c}|_{\Pi_{\sharp}}
$$
is regular in the ring ${\cal O}_{o,\Pi_{\sharp}}$.\vspace{0.1cm}

If the singularity $o\in F$ comes from a singularity of the original hypersurface $V$, then the condition (${\rm R}^*$ 2.1) follows from the condition (R2.1). If the singularity $o\in F$ comes from a non-singular point of the hypersurface $V$ (that is, $F$ is a section of $V$ by the hyperplane which is tangent to $V$ at this point), then (${\rm
R}^*$ 2.1) follows from (R1.1). In any case the condition (${\rm R}^*$
2.1) holds.\vspace{0.1cm}

(${\rm R}^*$ 2.2) The quadratic form $q^*_2$ is of rank at least
8.\vspace{0.1cm}

In our case this rank is at least 11.\vspace{0.1cm}

(${\rm R}^*$ 2.3) Considering $(u_1, \dots, u_M)$ as homogeneous coordinates $(u_1: \cdots : u_M)$ on ${\mathbb P}^{M-1}$, and the quadric hypersurface $\{q^*_2=0\}\subset{\mathbb P}^{M-1}$, let us construct the divisor $\{q^*_3|_{\{q^*_2=0\}}=0\}$. This divisor should not be a sum of three (not necessarily distinct) hyperplane sections of this quadric, taken from one  linear pencil.\vspace{0.1cm}

This condition follows from (R1.3), if the point $o\in F$ comes from a non-singular point of the hypersurface $V$, and from (R2.3), if $o\in F$ comes from a singular point of $V$.\vspace{0.1cm}

Thus we have checked that every hyperplane section of $V$ satisfies the regularity conditions of the paper \cite{Pukh15a}. Therefore, the global log canonical threshold of every fibre $F$ of the fibre space $\pi\colon V^+\to{\mathbb P}^1$ is equal to 1. The proof of Proposition 1.3 is complete.


\section{Regular hypersurfaces}

In this section we prove Theorem 2 and its immediate geometric implications. In Subsection 2.1 we consider all regularity conditions but the last one
(R2.3). In Subsection 2.2 we estimate the codimension for the violation of the condition (R2.3), using a technical fact which is shown in 2.3. This completes the proof of Theorem 2. In Subsection 2.4 we prove geometric facts which follow from the regularity conditions and will be used in the proof of Theorem 5 in Sections 3-5. For that purpose in Subsection 2.4 we briefly recall the technique of hypertangent divisors.\vspace{0.3cm}

{\bf 2.1. Violations of the regularity conditions.} Let us prove Theorem 2. We need to estimate the number of independent conditions which are imposed on the hypersurface $V$ (that is, on the coefficients of the polynomial $f$, defining this hypersurface) by violation of each of the six regularity conditions. Let us define the following polynomials of one real variable:
$$
\begin{array}{ccccl}
\gamma_{1.1}(t) & = & \gamma_{2.1}(t) & = & \frac12(t-8)(t-7)-14, \\
                &  & \gamma_{1.2}(t) & = & \frac12(t-11)(t-10)-10, \\
                &  & \gamma_{1.2}(t) & = & \frac12(t-2)(t^2-4t-27), \\
                &  & \gamma_{1.2}(t) & = & \frac12(t-11)(t-10)+1, \\
                &  & \gamma_{1.2}(t) & = & \frac12(t-9)(t-8)-11.
\end{array}
$$
Let ${\cal F}_{\alpha.\beta}\subset{\cal F}$ be the closed subset of hypersurfaces that do not satisfy the condition (R$\alpha.\beta$) at at least one point, where $\alpha\in\{1,2\}$ and $\beta\in\{1,2,3\}$. The following claim is true.\vspace{0.1cm}

{\bf Proposition 2.1.} {\it The following inequality holds:}
\begin{equation}\label{25.10.2017.2}
\mathop{\rm codim}({\cal F}_{\alpha.\beta}\subset{\cal F})\geqslant \gamma_{\alpha.\beta}(M).
\end{equation}

It is easy to see that for $M\geqslant 16$ the minimum of the values
$\gamma_{\alpha.\beta}(M)$ is attained for $\alpha=1$, $\beta=2$, which implies Theorem 2.\vspace{0.1cm}

{\bf Proof of Proposition 2.1.} Let us consider each of the conditions
(R$\alpha.\beta$) separately.\vspace{0.1cm}

{\bf The case $(\alpha,\beta)=(1,1)$.} The inequality (\ref{25.10.2017.2}) for these values of $\alpha$, $\beta$ is shown in
\cite[Proposition 2.5, Corollary 2.1]{Pukh16a} (p. 733-735).\vspace{0.1cm}

{\bf The case $(\alpha,\beta)=(1,2)$.} This is an elementary exercise.\vspace{0.1cm}

{\bf The case $(\alpha,\beta)=(1,3)$.} The inequality (\ref{25.10.2017.2}) for this case is shown in \cite[Section 2.6]{Pukh16a} (p. 735-736).\vspace{0.1cm}

{\bf The case $(\alpha,\beta)=(2,1)$.} Exactly the same arguments that prove the inequality (\ref{25.10.2017.2}) for $(\alpha,\beta)=(1,1)$, yield this inequality for $(\alpha,\beta)=(2,1)$, either, and the estimates turn out to be stronger, since the polynomials $q_2, q_3,\dots$ are restricted onto the linear space of a higher dimension.\vspace{0.1cm}

{\bf The case $(\alpha,\beta)=(2,2)$.} This is an elementary exercise.\vspace{0.1cm}

The only case which is non-trivial and was not considered in the previous papers, is the case $(\alpha,\beta)=(2,3)$. For that case, we give a complete detailed proof. We may (and will) assume the condition (R2.2) to be satisfied.\vspace{0.1cm}

For the convenience of our arguments set
$$
\gamma^*_{2.3}(t)=\frac12(t-7)(t-6)+1.
$$
Fix a point $o$ and a linear subspace
$\Pi\subset {\mathbb C}^{M+1}$ of codimension 2. Let us consider the closed set ${\cal F}_{2.3}(o,\Pi)\subset{\cal F}$, consisting of hypersurfaces $V$, such that

\begin{itemize}

\item contain the point  $o$,

\item are singular at that point,

\item do not satisfy the condition (R2.3) at that point for the subspace
$\Pi$.

\end{itemize}

Now we reduce the global statement (\ref{25.10.2017.2}) to the corresponding local statement.\vspace{0.1cm}

{\bf Proposition 2.2.} {\it The following inequality holds:}
$$
\mathop{\rm codim} ({\cal F}_{2.3}(o,\Pi)\subset{\cal F})\geqslant \gamma^*_{2.3}(M).
$$

Taking into account that the point $o$ runs through ${\mathbb P}$, and $\Pi$ varies in the $2(M-1)$-dimensional Grassmanian, and that the point $o$ lies on $V$ and is a singular point of that hypersurface, by an elementary dimension count we check that Proposition 2.2 implies the inequality (\ref{25.10.2017.2}) for $(\alpha,\beta)=(2,3)$. Q.E.D. for Proposition 2.1.\vspace{0.3cm}


{\bf 2.2. Violation of the condition (R2.3).} Let us prove Proposition 2.2. The symbol ${\cal P}_{k,\Pi}$ stands for the space of homogeneous polynomials of degree $k$ on ${\mathbb P}(\Pi)\cong{\mathbb P}^{M-2}$. The restrictions of the polynomials $q_2$ and $q_3$ onto ${\mathbb P}(\Pi)$ are denoted by the symbols ${\bar q}_2$ and ${\bar q}_3$, respectively, and the set of their common zeros $\{{\bar q}_2={\bar q}_3=0\}\subset{\mathbb P}(\Pi)$ by the symbol $Z({\bar q}_2,{\bar q}_3)$. It is obvious that the codimension of the set ${\cal F}_{2.3}(o,\Pi)$, which we need to estimate, is equal to the codimension of the set ${\cal B}\subset{\cal P}_{2,\Pi}\times{\cal P}_{3,\Pi}$ of pairs $({\bar q}_2,{\bar q}_3)$, such that the set $Z({\bar q}_2,{\bar q}_3)$ is not an irreducible reduced factorial complete intersection of type $2\cdot 3$ in ${\mathbb P}(\Pi)$.\vspace{0.1cm}

We note at once that the quadratic form ${\bar q}_2$ is by the condition
(R2.2) of rank at least 9, so that the quadric
$\{{\bar q}_2=0\}\subset{\mathbb P}(\Pi)$ is for sure factorial. It is easy to compute that for a fixed form ${\bar q}_2$ of rank $\geqslant 9$ the set of cubic polynomials ${\bar q}_3\in{\cal P}_{3,\Pi}$, such that the divisor
$\{{\bar q}_3|_{\{{\bar q}_2=0\}}=0\}$ on the quadric $\{{\bar q}_2=0\}$ is non-reduced or reducible has codimension
$\frac16 M(M+1)(M-4)$ in ${\cal P}_{3,\Pi}$. Since this is much higher than $\gamma_{2.3}(M)$, we may (and will) assume that the set
$Z({\bar q}_2,{\bar q}_3)$ is irreducible, reduced and of codimension 2 in
${\mathbb P}(\Pi)$. It remains to consider the condition for this set to be  factorial. Let $p\in Z({\bar q}_2,{\bar q}_3)$ be an arbitrary point and
$(u_1, \dots, u_{M-2})$ a system of affine coordinates on ${\mathbb P}(\Pi)$ with the origin at the point $p$. Let ${\mathbb P}(\Pi)^+\to{\mathbb P}(\Pi)$ be the blow up of the point $p$ with the exceptional divisor $E_p\cong{\mathbb P}^{M-3}$, equipped with the natural homogeneous coordinates $(u_1: \cdots : u_{M-2})$. The affine polynomials in the (non-homogeneous) variables $u_*$, corresponding to to the homogeneous polynomials ${\bar q}_2,{\bar q}_3$, we denote, somewhat abusing the notations, by the same symbols ${\bar q}_2,{\bar q}_3$. We get
$$
\begin{array}{ccccccc}
{\bar q}_2 & = & {\bar q}_{2,1} & + & {\bar q}_{2,2}, &  & \\
& & & & & & \\
{\bar q}_3 & = & {\bar q}_{3,1} & + & {\bar q}_{3,2} & + & {\bar q}_{3,3},
\end{array}
$$
where ${\bar q}_{i,j}$ are homogeneous of degree $j$. We say that the point $p$ is a {\it correct bi-quadratic point} of the set $Z({\bar q}_2,{\bar q}_3)$, if ${\bar q}_{2,1}\equiv {\bar q}_{3,1}\equiv 0$, and the closed set
$$
\{{\bar q}_{2,2}={\bar q}_{3,2}=0\}\subset E_p
$$
is an irreducible reduced complete intersection of codimension 2 in
$E_p\cong{\mathbb P}^{M-3}$, and moreover,
$$
\mathop{\rm dim} \mathop{\rm Sing}\{{\bar q}_{2,2}={\bar q}_{3,2}=0\}\leqslant M-9.
$$
Let ${\cal X}\subset {\cal P}_{2,\Pi}\times{\cal P}_{3,\Pi}$ be the set of pairs such that $Z({\bar q}_2,{\bar q}_3)$ is irreducible, reduced, of codimension 2 in ${\mathbb P}(\Pi)$, and its every point
$p\in Z({\bar q}_2,{\bar q}_3)$

\begin{itemize}

\item either is non-singular,

\item or is a quadratic singularity of rank $\geqslant 5$,

\item or is a correct bi-quadratic point.

\end{itemize}

By Grothendieck's theorem \cite{CL} for the pair
$({\bar q}_2,{\bar q}_3)\in{\cal X}$ the complete intersection
$Z({\bar q}_2,{\bar q}_3)$ is factorial, so that
${\cal B}\cap{\cal X}=\emptyset$ and in order to prove Proposition 2.2, it is sufficient to show that the codimension of the complement to the set ${\cal X}$ in ${\cal P}_{2,\Pi}\times {\cal P}_{3,\Pi}$ is at least
$\gamma^*_{2.3}(M)$.\vspace{0.1cm}

Recall that $\mathop{\rm rk} {\bar q}_2\geqslant 9$, and the complete intersection $Z({\bar q}_2,{\bar q}_3)$ is irreducible and reduced. Fix a point $p\in Z({\bar q}_2,{\bar q}_3)$. Assume first that at least one of the linear forms ${\bar q}_{2,1}$, ${\bar q}_{3,1}$ is not identically zero, but these forms are linearly dependent. If, moreover, ${\bar q}_{2,1}\equiv 0$, then ${\bar q}_{3,1}\not\equiv 0$ and the point $p$ is a quadratic singularity of rank
$$
\mathop{\rm rk} {\bar q}_{2,2}|_{\{{\bar q}_{3,1}=0\}}\geqslant 9-4=5,
$$
which is what we need. If ${\bar q}_{2,1}\not\equiv 0$, then there is a unique constant $\lambda$, such that
${\bar q}_{3,1}=\lambda {\bar q}_{2,1}$. In that case $p$ is a quadratic singularity of rank
$$
\mathop{\rm rk} ({\bar q}_{3,2}-\lambda{\bar q}_{2,2})|_{\{{\bar q}_{2,1}=0\}}.
$$
If this rank $\leqslant 4$, then for a fixed polynomial ${\bar q}_2$ we get $\frac12 (M-7)(M-6)$ independent conditions for the polynomial
${\bar q}_3$. Taking into account that the constant $\lambda$ varies in a 1-dimensional family, there is the dependence ${\bar q}_{3,1}=\lambda{\bar q}_{2,1}$, the point $p$ varies in ${\mathbb P}^{M-2}$ and the polynomials
${\bar q}_2$, ${\bar q}_3$ vanish at that point, we get precisely the codimension $\gamma^*_{2.3}(M)$ for the violation of the condition about the rank of quadratic points.\vspace{0.1cm}

It remains to consider the case ${\bar q}_{2.1}\equiv {\bar q}_{3.1}\equiv 0$ and estimate the codimension for the violation of the condition about bi-quadratic points. For this purpose we state and solve the following general problem. Let ${\cal P}_{2,N+1}$ be the space of quadratic forms on ${\mathbb P}^N$, where $N\geqslant 8$. Let
${\cal Y}\subset{\cal P}^{\times 2}_{2,N+1}$ be the set of pairs
$(g_1,g_2)$, such that the closed set of common zeros
$$
Z(g_1,g_2)=\{g_1=g_2=0\}\subset{\mathbb P}^N
$$
is an irreducible reduced complete intersection of codimension 2, and moreover,
$$
\mathop{\rm codim}(\mathop{\rm Sing} Z(g_1,g_2)\subset{\mathbb P}^N)\geqslant 6.
$$
The following fact is true.\vspace{0.1cm}

{\bf Proposition 2.3.} {\it The codimension of the complement to ${\cal Y}$ in ${\cal P}^{\times 2}_{2,N+1}$ is at least} $\frac12 (N-4)(N-3)-2$.\vspace{0.1cm}

Let us complete the proof of Proposition 2.2. Setting in Proposition 2.3
$N=M-3$, we get that violation of the condition about bi-quadratic points at the fixed point $p$ gives the codimension $\frac12 (M-7)(M-6)-2$. Now the standard dimension count (taking into account that
${\bar q}_{2,1}\equiv{\bar q}_{3,1}\equiv 0$, and also the conditions
$p\in Z({\bar q}_2,{\bar q}_3)$ and the variation of the point $p$) completes the proof of Proposition 2.2.\vspace{0.3cm}


{\bf 2.3. Complete intersections of two quadrics.} Let us prove Proposition 2.3. Since the set of quadratic forms of rank $\leqslant 4$ has codimension $\frac12 (N-3)(N-2)$ in ${\cal P}_{2,N+1}$, we may assume that
$\mathop{\rm rk} g_i\geqslant 5$ for $i=1,2$, so that the quadric
$\{g_1=0\}$ is factorial. If the condition of irreducibility and reducedness of the set $Z(g_1,g_2)$ is violated, this imposes on $g_2$ a lot more conditions than the required $\frac12 (N-4)(N-3)-2$. Thus only the condition about the singularities of the set $Z(g_1,g_@)$ needs to be considered. We argue as in \cite[Section 3.3]{EvP2017}, somewhat improving the estimate obtained in that paper. The key observation (used in \cite{EvP2017}) is that if $p\in\mathop{\rm Sing} Z(g_1,g_2)$, then for some
$\lambda_1,\lambda_2$ (where $(\lambda_1,\lambda_2)\neq (0,0)$) the point $p$ is a singular point of the quadric $\{\lambda_1g_1+\lambda_2g_2=0\}$. In order to obtain a somewhat more precise, than in \cite{EvP2017}, estimate for the codimension of the set of ``incorrect'' pairs, we have to consider several cases. For a quadratic polynomial $g\in{\cal P}_{2,N+1}$ the symbol
$C(g,\leqslant k)$ stands for the cone with the vertex $g$, the base of which is the set of quadratic forms of rank $\leqslant k$. The vertex can lie on the base: $C(g,\leqslant k)$ is the closure
$$
\overline{\mathop{\bigcup}\limits_{\mathop{\rm rk} h\leqslant k}
[g,h]},
$$
where $[g,h]=\{\lambda g+\mu h\,|\, \lambda,\mu\in{\mathbb C}\}$. Obviously,
$$
\mathop{\rm codim}(C(g,\leqslant k)\subset{\cal P}_{2,N+1})\geqslant
\frac12 (N-k+1)(N-k+2)-2.
$$

{\bf Case 1: $\mathop{\rm rk} g_1=5$.} We have $\frac12 (N-4)(N-3)$ independent conditions for $g_1$. If for a fixed quadric $g_1$ we have, into the bargain, $g_2\in C(g_1,\leqslant 6)$, this gives in addition
$\frac12 (N-5)(N-4)-2$ independent conditions for $g_2$, and we get the total $(N-4)^2-2$ independent conditions for the pair $(g_1,g_2)$, which is much higher than what we need. Therefore we may assume that $g_2\not\in C(g_1,\leqslant 6)$.\vspace{0.1cm}

This implies that in the pencil $\{\lambda_1 g_1+\lambda_2 g_2=0\}$ all quadrics, apart from $g_1$, are of rank $\geqslant 7$. Therefore, the codimension of the set
$$
\overline{\mathop{\bigcup}\limits_{(\lambda_1:\lambda_2)\neq (1:0)}
\mathop{\rm Sing} \{\lambda_1 g_1+\lambda_2 g_2=0\}}
$$
in ${\mathbb P}^N$ is at least 6. On the other hand,
$\mathop{\rm Sing} \{g_1=0\}$ is a $(N-5)$-dimensional subspace in
${\mathbb P}^N$, so that the condition
$$
g_2|_{\mathop{\rm Sing} \{g_1=0\}}\equiv 0
$$
gives for $g_2$ the codimension $\frac12 (N-4)(N-3)$ (for a fixed
$g_1$), and for the pair $(g_1,g_2)$ the codimension $(N-4)(N-3)$. Removing this set of high codimension, we may assume that
$$
g_2|_{\mathop{\rm Sing} \{g_1=0\}}\not\equiv 0,
$$
but then the set
$$
\mathop{\rm Sing} Z(g_1,g_2)\cap \mathop{\rm Sing} \{g_1=0\}
$$
is of dimension at most $N-6$. We get finally that the codimension of the set $\mathop{\rm Sing} Z(g_1,g_2)$ is at least 5 in
${\mathbb P}^N$ in the case under consideration (for the pairs
$(g_1,g_2)\in{\cal P}^{\times 2}_{2,N+1}$, lying outside a closed subset of high codimension).\vspace{0.1cm}

{\bf Case 2: $\mathop{\rm rk} g_1=6$.} Here we argue word for word as in Case 1. The only difference is that we get somewhat fewer
($\frac12(N-5)(N-4)$) independent conditions for $g_1$. Together with the conditions for the form $g_2$ we get the total codimension of the set of pairs $(g_1,g_2)$ to be much higher than we need. In this case we do not need to exclude the option
$g_2|_{\mathop{\rm Sing} \{g_1=0\}}\equiv 0$.\vspace{0.1cm}

{\bf Case 3: $\mathop{\rm rk} g_1\geqslant 7$.} This is the case of general position for $g_1$, so that no conditions are imposed on $g_1$. For a fixed form $g_1$ the condition
$$
g_2\in C(g_1,\leqslant 5)
$$
imposes on $g_2$ precisely $\frac12 (N-4)(N-3)-2$ independent conditions. This is the estimate that we need. If $g_2\not\in C(g_1,\leqslant 5)$, then a general quadric in the pencil $\{\lambda_1 g_1 +\lambda_2 g_2=0\}$ is of rank $\geqslant 7$ and at most finitely many quadrics are of rank 6. This implies that
$$
\mathop{\rm codim}(\mathop{\rm Sing}(g_1,g_2)\subset{\mathbb P}^N)\geqslant 6,
$$
which is what we need.\vspace{0.1cm}

Proof of Proposition 2.3 is complete. Q.E.D.\vspace{0.3cm}


{\bf 2.4. Hypertangent divisors at a singular point.} The technique of hypertangent divisors makes it possible to obtain very strong upper bounds for the multiplicity of irreducible subvarieties (and thus of effective cycles) of the hypersurface $V$ at a given point. For the full details of this technique for various types of Fano varieties see \cite[Chapter 3]{Pukh13a}. For a non-singular point $o\in V$ see \cite[Sections 2.5 and 4.1]{Pukh16a}. Here we briefly consider the case of a singular point $o\in V$ and give the estimates that will be used later.\vspace{0.1cm}

In the notations of Subsection 0.2, let $\Pi\subset{\mathbb C}^{M+1}$ be a linear subspace of codimension $c\in\{1,2,3\}$, $\overline{\Pi}\subset {\mathbb P}$ its closure and $V_{\Pi}=V\cap\overline{\Pi}$ the corresponding section of the hypersurface $V$. For an irreducible subvariety $Y\subset V$ the symbol
$$
\frac{\mathop{\rm mult}\nolimits_o}{\mathop{\rm deg}} Y
$$
stands for the ratio
${\mathop{\rm mult}\nolimits_o Y}/{\mathop{\rm deg}} Y$, where the degree is understood in the usual sense as the degree in the projective space
${\mathbb P}$. A similar symbol will be used in the sequel for an effective equidimensional cycle, either.\vspace{0.1cm}

{\bf Proposition 2.4.} {\it For an irreducible subvariety
$Y\subset V_{\Pi}$ of codimension $a$, where $1\leqslant a\leqslant M-c-1$, the following inequality holds:}
\begin{equation}\label{04.11.2017.1}
\frac{\mathop{\rm mult}\nolimits_o}{\mathop{\rm deg}} Y\leqslant \frac{a+2}{M+1-c}.
\end{equation}

{\bf Remark 2.1.} (i) In the case $c=0$, when $V_{\Pi}=V$, the inequality
(\ref{04.11.2017.1}) for $c=1$ implies the estimate
$$
\frac{\mathop{\rm mult}\nolimits_o}{\mathop{\rm deg}} Y\leqslant\frac{a+2}{M}
$$
for any irreducible subvariety $Y\subset V$ of codimension
$a$, where $1\leqslant a \leqslant M-2$.\vspace{0.1cm}

(ii) If $Y$ is a curve, that is, $a=M-c-1$, then the inequality
(\ref{04.11.2017.1}) is trivial:
$\mathop{\rm mult}\nolimits_o Y\leqslant \mathop{\rm deg} Y$.\vspace{0.1cm}

{\bf Proof of Proposition 2.4.} We give a sketch of the arguments, which are absolutely standard, see \cite[Chapter 3]{Pukh13a}. In the notations of Subsection 0.2 set for $i=2,\dots, M-1$
$$
f_i=q_2+\cdots +q_i
$$
and for every $j=2,\dots, M-1$ construct the linear system
$$
\Lambda_j=\left\{\left(\left.\sum^j_{\alpha=2} f_{\alpha} s_{j-\alpha}\right)\right|_{V_{\Pi}}=0\right\},
$$
where $s_{\beta}$ are homogeneous polynomials of degree
$\beta\geqslant 0$ in the variables $z_*$, which run through the spaces ${\cal P}_{\beta,M+1}$ independently of each other. The linear system $\Lambda_j$ is called the $j$-the hypertangent linear system on $V_{\Pi}$ at the point $o$. The regularity condition (R2.1) implies that for
$j\in\{2,3,\dots, M-c\}$ the equality
$$
\mathop{\rm codim}\nolimits_o(\mathop{\rm Bs} \Lambda_j\subset V_{\Pi})=j-1
$$
holds, where the symbol $\mathop{\rm codim}\nolimits_o$ stands for the codimension in a neighborhood of the point $o$. If $Y\not\ni o$, then the inequality (\ref{04.11.2017.1}) is trivial. Assume that $Y\ni o$. In that case for a general hypertangent divisor $T_{a+2}\in\Lambda_{a+2}$ we get:
$$
Y\not\subset |T_{a+2}|
$$
(the vertical lines mean the support of the divisor), so that the effective cycle of the scheme-theoretic intersection $(Y\circ T_{a+2})$ of codimension $(a+1)$ on $V_{\Pi}$ is well defined. Since
$\Lambda_j\subset |jH_{\Pi}|$ (where $H_{\Pi}$ is the class of a hyperplane section of $V_{\Pi}$) and $\mathop{\rm mult}\nolimits_o\Lambda_j\geqslant (j+1)$, we get
$$
\frac{\mathop{\rm mult}\nolimits_o}{\mathop{\rm deg}} (Y\circ T_{a+2})\geqslant
\frac{a+3}{a+2}\cdot \frac{\mathop{\rm mult}\nolimits_o}{\mathop{\rm deg}} Y,
$$
whence the inequality (\ref{04.11.2017.1}) is obtained by decreasing induction on $a$ (see Remark 2.1, (ii)). For the details, see
\cite[Chapter 3]{Pukh13a}. Q.E.D. for Proposition 2.4.\vspace{0.1cm}

The hypertangent system $\Lambda_2$ is not mobile: it consists of the unique divisor $T_2=\{q_2|_{V_{\Pi}}=0\}$. By the condition (R2.2) and the factoriality of $V_{\Pi}$, this divisor is irreducible. For that reason, the claim of Proposition 2.4 for divisors can be made slightly more precise.\vspace{0.1cm}

{\bf Proposition 2.5.} {\it Let $Y\subset V_{\Pi}$ be a prime divisor, and moreover, $Y\neq T_2$. Then the following inequality holds:}
\begin{equation}\label{04.11.2017.2}
\frac{\mathop{\rm mult}\nolimits_o}{\mathop{\rm deg}} Y\leqslant \frac{8}{3(M+1-c)}.
\end{equation}

{\bf Proof.} Let us apply Proposition 2.4 to the effective cycle of the intersection $(Y\circ T_2)$ of codimension 2 on $V_{\Pi}$. Q.E.D.


\section{Maximal singularities of the system $\Sigma$}

In this section we begin to study the maximal singularity
$E^*$. First (Subsection 3.1) we show that the centre $B^*$ of this singularity is contained in $\mathop{\rm Sing} V$. In order to do this, we check that the arguments of \cite{Pukh16a} exclude a maximal singularity, the centre of which is not contained in $\mathop{\rm Sing} V$. After that (Subsection 3.2), using inversion of adjunction, we derive from the existence of the singularity $E^*$ the existence of a, generally speaking, another singularity $E$ of the linear system $\Sigma$, the centre of which is a point $o\in\mathop{\rm Sing} V$, and moreover, $E$ has some good properties (which may not be satisfied for $E^*$). Then we classify the types of the singularity $E$. Finally, in Subsection 3.3 we recall the technique of counting multiplicities which makes use of combinatorial invariants of the oriented graph associated with the singularity $E$.\vspace{0.3cm}

{\bf 3.1. The linear system $\Sigma$ at non-singular points of the hypersurface $V$.} Starting from this moment, we assume that the mobile linear system $\Sigma\subset |2H|$ satisfies the inequality
$$
\mathop{\rm mult}\nolimits_{P\cap V}\Sigma\leqslant n
$$
for every linear subspace $P\subset{\mathbb P}$ of codimension 2. On the other hand, the system $\Sigma$ has the maximal singularity $E^*$ (see Subsection 0.3). Theorem 5 will be shown if we derive a contradiction from this fact. We will do it in several steps, excluding the possible types of the maximal singularity $E^*$ one after another. The first step is given by the following statement.\vspace{0.1cm}

{\bf Proposition 3.1.} {\it The centre $B^*$ of the maximal singularity $E^*$ on $V$ is contained in the singular locus}
$\mathop{\rm Sing} V$.\vspace{0.1cm}

{\bf Proof.} Let us assume the converse:
$B^*\not\subset\mathop{\rm Sing} V$ and show that the arguments of
\cite{Pukh16a} exclude this option. Let us consider separately the three cases:\vspace{0.1cm}

{\bf Case 1.} $\mathop{\rm codim} (B^*\subset V)=2$.\vspace{0.1cm}

{\bf Case 2.} $\mathop{\rm codim} (B^*\subset V)\in\{3,4,\dots, 9\}$.\vspace{0.1cm}

{\bf Case 3.} $\mathop{\rm codim} (B^*\subset V)\geqslant 10$.\vspace{0.1cm}

The maximal singularity is excluded in each of these three cases in a different way. Note that the inequality
$\mathop{\rm mult}\nolimits_{B^*}\Sigma>n$ holds.\vspace{0.1cm}

Assume that {\bf Case 1} takes place. In this case $B^*$ is a maximal subvariety of the system $\Sigma$ and, arguing in a word for word the same way as in \cite[Section 3.1]{Pukh16a}, we conclude that
$\langle B^*\rangle={\mathbb P}$ (since every hyperplane section of the hypersurface $V$ is a factorial variety). Furthermore, in the notations of \cite[Section 3.1]{Pukh16a} we conclude that
$\mathop{\rm Sec}(B^*)=\pi_{{\mathbb P}}(\mathop{\rm\bf Sec}(B^*))={\mathbb P}$, that is, the claim of \cite[Proposition 3.1]{Pukh16a} is true in our case, either.\vspace{0.1cm}

Indeed, the proof of that claim, given in
\cite[Section 3.3]{Pukh16a}, makes use of only one fact, that
$B^*$ is contained in a non-singular hypersurface, which does not contain cones over a positive-dimensional base (that is, cones of dimension
$\geqslant 2$). In our case the inequality
$$
\mathop{\rm codim}(\mathop{\rm Sing} V\subset V)\geqslant 12
$$
holds, which follows from the regularity condition (R2.2), so that for a general linear subspace $\Pi\subset{\mathbb P}$ of dimension 12 the hypersurface $V_{\Pi}=V\cap \Pi$ is non-singular and $B^*\cap\Pi$ is an irreducible subvariety of codimension 2 on $V_{\Pi}$. Furthermore, $V_{\Pi}$ does not contain cones of dimension $\geqslant 2$, because $V$ does not contain any (by the conditions (R1.1) and (R2.1), there are at most finitely many lines on $V$ through every point of $V$). Therefore, we have
$\mathop{\rm Sec}(B^*\cap\Pi)=\Pi$, which implies that
$\mathop{\rm Sec}(B^*)={\mathbb P}$, as we claimed.\vspace{0.1cm}

Now the arguments of \cite[Sections 3.1,3.2]{Pukh16a}, excluding a maximal subvariety of codimension 2 --- in our case it is denoted by $B^*$ --- work word for word. The curves $C_{\pm}$, $R$ do not touch the set
$\mathop{\rm Sing} V$ for a general point $x\in{\mathbb P}$ and a general curve $\Gamma$, see \cite[Section 3.1]{Pukh16a} (p. 739-740). Therefore, Case 1 can not take place.\vspace{0.1cm}

Assume that {\bf Case 2} takes place. Here we argue as in the proof of \cite[Lemma 4.1]{Pukh16a}. Let $Z=(D_1\circ D_2$ be the self-intersection of the linear system $\Sigma$, where $D_1,D_2\in\Sigma$ are general divisors. The $4n^2$-inequality holds:
$$
\mathop{\rm mult}\nolimits_{B^*}Z>4n^2.
$$
Again let $\Pi$ be a general linear subspace of dimension 12. Then
$V_{\Pi}=V\cap \Pi\subset\Pi\cong{\mathbb P}^{12}$ is a non-singular hypersurface and $B\cap \Pi$ an irreducible subvariety,
$\mathop{\rm dim} B^*\cap \Pi\geqslant 2$, and the effective cycle
$Z_{\Pi}=(Z\circ V_{\Pi})$ of codimension 2 on $V_{\Pi}$ satisfies the inequality
$$
\mathop{\rm mult}\nolimits_{B^*\cap\Pi}Z_{\Pi}>4n^2.
$$
However $Z\sim 4n^2H^2$, so that $Z_{\Pi}\sim 4n^2H^2_{\Pi}$, where
$H_{\Pi}$ is the class of a hyperplane section of the hypersurface $V_{\Pi}$. Now as in the proof of \cite[Lemma 4.1]{Pukh16a} the reference to
\cite[Proposition 5]{Pukh02f} gives a contradiction, excluding Case 2.\vspace{0.1cm}

Assume that {\bf Case 3} takes place. Let us check that the arguments of
\cite[Sections 4-6]{Pukh16a} exclude this case. Let $o\in B^*$ be a point of general position, in particular, $o\not\in\mathop{\rm Sing} V$. The claims of \cite[Propositions 4.1-4.3]{Pukh16a} are true: they follow from the regularity condition which is identical to the condition (R1.1). Furthermore, the $8n^2$-inequality (\cite[Proposition 4.4]{Pukh16a}) is a local fact and for that reason is true in our case. Indeed, we may assume that $B^*$ has maximal dimension among all centres of maximal singularities, containing the point $o$ (since Cases 1 and 2 are already excluded). Since for a general linear subspace $R\subset{\mathbb P}$ of dimension 11, containing the point $o$, this point is an {\it isolated} centre of a non-canonical singularity of the pair $(V_R,\frac{1}{n}\Sigma_R)$, where $V_R=V\cap R$ and
$\Sigma_R=\Sigma|_{V_R}\subset |2nH_R|$ is a mobile linear system
($H_R$ means the class of a hyperplane section of the hypersurface
$V_R$). Therefore, for a general linear subspace $R_1$,
$o\in R_1\subset R$, where $5\leqslant \mathop{\rm dim} R_1\leqslant 10$, by inversion of adjunction the pair
$$
(V\cap R_1, \frac{1}{n}\Sigma|_{R_1})
$$
is not {\it log} canonical at the point $o$. This implies the
$8n^2$-inequality (see \cite[Chapter 2, Section 4]{Pukh13a}) and the existence of a hyperplane section $P\ni o$, such that

\begin{itemize}

\item the linear system $\Sigma_P=\Sigma|_P\subset |2nH_P|$ is mobile,
($H_P$ is the class of a hyperplane section of the variety $P$),

\item its self-intersection $Z_P=(Z\circ P)$ satisfies the inequality
$\mathop{\rm mult}\nolimits_o Z_P> 8n^2$,

\end{itemize}
see \cite[p. 753]{Pukh16a}. Furthermore, by what was said above, the claim of \cite[Proposition 4.5]{Pukh16a} is true (the proof works word for word) and the subsequent arguments of \cite[Section 4.2]{Pukh16a}, showing the existence of non log canonical singularities of the pairs
$\Box$ and $\Box^*$ (see \cite[p. 755]{Pukh16a}). The remaining part of the paper \cite{Pukh16a} (starting from Section 4.3 and up to the end) excludes these singularities of the pairs $\Box$ and $\Box^*$, and the arguments of \cite{Pukh16a} work word for word, without any modifications, except for the only place: \cite[Section 5.4]{Pukh16a}, where formally speaking the regularity condition (R2) is used, which requires that the rank of the quadratic form
$$
q_2|_{\{q_1=0\}}
$$
is at least $M-[\frac12 (\sqrt{8M+1}-1)]$. This condition is stronger than the regularity condition (R1.2) used in this paper. However, that condition
(R2) is used in \cite[Section 5.4]{Pukh16a} only once --- in the proof of Corollary 5.1, and it is easy to see from the proof that the condition (R2) is unnecessarily strong: the inequality
$$
M-10>M-\mathop{\rm rk}q_2|_{\{q_1=0\}},
$$
which is equivalent to the condition (R1.2) of the present paper, is sufficient. Thus all arguments of the paper \cite{Pukh16a}, excluding the maximal singularity in Sections 4-6, work without modifications in our Case 3 and exclude the maximal singularity, the centre of which is not contained in $\mathop{\rm Sing} V$. This completes the proof of Proposition 3.1.\vspace{0.3cm}


{\bf 3.2. The linear system $\Sigma$ at singular points of the hypersurface
$V$.} Let us fix a maximal singularity $E^*$, the centre $B^*$ of which has the maximal dimension among all centres of maximal singularities of the linear system $\Sigma$. We have $B^*\subset\mathop{\rm Sing} V$. Let
$o\in B^*$ be a point of general position. For a general 13-dimensional linear subspace $\Pi\subset{\mathbb P}$, where $o\in\Pi$, the pair
$$
(V_{\Pi},\frac{1}{n}\Sigma_{\Pi}),
$$
where $V_{\Pi}=V\cap \Pi$ and $\Sigma_{\Pi}$ is the restriction of $\Sigma$ onto $V_{\Pi}$, has the point $o$ as an isolated centre of a non canonical singularity, that is, this pair is canonical outside the point $o$ in a neighborhood of that point. By inversion of adjunction for a general proper subspace $\Pi_1\subset\Pi$, containing the point $o$, the pair
$$
(V\cap\Pi_1,\frac{1}{n}\Sigma|_{V\cap \Pi_1})
$$
is not log canonical at the point $o$, but canonical outside that point in a neighborhood of that point.\vspace{0.1cm}

Let $\varphi_{{\mathbb P}}\colon{\mathbb P}^+\to{\mathbb P}$ be the blow up of the point $o$ with the exceptional divisor $E_{\mathbb P}\cong{\mathbb P}^M$, and $\varphi\colon V^+\to V$ the restriction of that blow up onto the hypersurface $V$. The exceptional divisor $Q=V^+\cap E_{\mathbb P}$ of the blow up $\varphi$ is by the condition (R2.2) a quadric hypersurface of rank at least 13, embedded in $E_{\mathbb P}$. The symbol $H_Q$ stands for the class of a hyperplane section of $Q$. Every irreducible subvariety
$R\subset Q$ of codimension $\leqslant 5$ is numerically equivalent to the class
$$
\mathop{\rm d}\nolimits_Q (R)H_Q^{\mathop{\rm codim} (R\subset Q)}
$$
for some $\mathop{\rm d}\nolimits_Q (R)\in{\mathbb Z}_+$; by linearity, the integer-valued function $\mathop{\rm d}\nolimits_Q(\cdot)$ extends for all equidimensional cycles of codimension $\leqslant 5$.\vspace{0.1cm}

Furthermore, let $R\subset V$ be an irreducible subvariety of codimension
$\leqslant 5$. We get the numerical equivalence
$$
R\sim \mathop{\rm d} (R)H^{\mathop{\rm codim} (R\subset V)}
$$
for some $\mathop{\rm d} (R)\in{\mathbb Z}_+$; again, this integer-valued function extends by linearity to all equidimensional cycles of this codimension. Let $R^+\subset V^+$ be the strict transform of $R$ on $V^+$. Obviously,
$$
R^+\sim \varphi^*R-\mathop{\rm m}(R) H_Q^{\mathop{\rm codim} (R\subset V)-1},
$$
where $H_Q^0=Q$ and
$$
\mathop{\rm m}(R)=\mathop{\rm d}\nolimits_Q (R^+\cap Q)
$$
again extends to equidimensional cycles. For simplicity of notations we will often omit the pull back symbol: for instance, we write $R$ instead of
$\varphi^* R$. The obvious equalities hold:
$$
\mathop{\rm deg} R=M\mathop{\rm d}(R),\quad \mathop{\rm mult}\nolimits_o R=2\mathop{\rm m}(R).
$$
For a general divisor $D\in\Sigma$ set
$$
\nu=\mathop{\rm m}(D),
$$
that is, $D^+\sim D-\nu Q$. Proposition 2.5 implies the inequality
$\nu\leqslant\frac83 n$. Consider the self-intersection
$Z=(D_1\circ D_2)$ of the mobile system $\Sigma$. By construction,
$\mathop{\rm d}(Z)=4n^2$. The singularity $o\in V$ satisfies the assumptions of the main theorem of \cite{Pukh2017a}, therefore the inequality
$$
\mathop{\rm m}(Z)>4n^2
$$
holds (which, unfortunately, is insufficient for the exclusion of the maximal singularity).\vspace{0.1cm}

Let $\Pi\subset{\mathbb P}$ be a general 6-dimensional subspace, containing the point $o$, $V_{\Pi}=V\cap\Pi$ and $\Sigma_{\Pi}=\Sigma|_{V_{\Pi}}$. The symbol $H_{\Pi}$ stands for the class of a hyperplane section of the hypersurface $V_{\Pi}$. Obviously, $V_{\Pi}$ has a unique singular point --- the non-degenerate quadratic point $o$. Let $V_{\Pi}^+$ be the strict transform of $V_{\Pi}$ on $V^+$, that is,
$$
\varphi_{\Pi}\colon V^+_{\Pi}\to V_{\Pi}
$$
is the blow up of the point $o$ with the exceptional divisor
$Q_{\Pi}=V_{\Pi}^+\cap E_{\mathbb P}$, which is a non-singular 4-dimensional quadric in $\Pi^+\cap E_{\mathbb P}\cong{\mathbb P}^5$. We get
$$
\Sigma_{\Pi}^+\subset |2nH_{\Pi}-\nu Q_{\Pi}|.
$$
Since for the discrepancy of the exceptional divisor $Q_{\Pi}$ we have the equality
$$
a(Q_{\Pi},V_{\Pi})=3
$$
and as we mentioned above, $\nu\leqslant 3n$ and the pair
$(V_{\Pi},\frac{1}{n}\Sigma_{\Pi})$ is not log canonical at the point
$o$, we obtain the following fact: the pair
$$
\Box_{\Pi}=
\left(V_{\Pi}^+,\frac{1}{n}\Sigma_{\Pi}^++\frac{\nu-3n}{n}Q_{\Pi}\right)
$$
is not log canonical, and moreover, the centre of every non log canonical singularity of the pair ${\Box}_{\Pi}$, intersecting the exceptional quadric $Q_{\Pi}$, is contained in it. Furthermore, the pair ${\Box}_{\Pi}$ satisfies the assumptions of the connectedness principle of Shokurov and
Koll\'{a}r \cite[Section 17.4]{Kol93} and for that reason the union
$\mathop{\rm LCS}({\Box}_{\Pi})$ of all centres of non log canonical singularities of the pair ${\Box}_{\Pi}$ is a connected closed subset of the quadric $Q_{\Pi}$. We say that {\it Case $\alpha$} takes place, where
$\alpha\in\{1,2,3,4\}$, if
$$
\mathop{\rm codim}(\mathop{\rm LCS}({\Box}_{\Pi})\subset Q_{\Pi})=\alpha,
$$
where if $\alpha=4$, then $\mathop{\rm LCS}({\Box}_{\Pi})$ is one point. For the original hypersurface $V$ this means that the pair
$$
\Box=\left(V^+,\frac{1}{n}\Sigma^++\frac{\nu-3n}{n}Q\right)
$$
has a non log canonical singularity $E$ (a certain exceptional divisor over
$V^+$), the centre $B$ of which is contained in the exceptional quadric
$Q$, and moreover $\mathop{\rm codim} (B\subset Q)=\alpha\in\{1,2,3,4\}$, and if $\alpha=4$, then $B\subset E_{\mathbb P}$ is a linear subspace of dimension $M-5$.\vspace{0.1cm}

{\bf Proposition 3.2.} {\it The case $\alpha=4$ is impossible.}\vspace{0.1cm}

{\bf Proof.} Indeed, the quadric of rank $\geqslant 13$ can not contain a linear subspace of codimension 4 (the codimension is meant with respect to this quadric). Q.E.D. for the proposition.\vspace{0.1cm}

{\bf Proposition 3.3.} {\it The case $\alpha=1$ is impossible.}\vspace{0.1cm}

{\bf Proof.} Assume the converse: $B\subset Q$ is a prime divisor. Arguing in precisely the same way as in the proof of
\cite[Proposition 4.6]{Pukh16a} (which, in its turn, repeats the arguments in \cite[Chapter 2, Section 4.1]{Pukh13a} and in
\cite{Ch06b}), and taking into account the fact that the point $o$ is a double point of the hypersurface $V$, we obtain the following inequality for the self-intersection $Z$:
$$
\mathop{\rm mult}\nolimits_o Z>2\nu^2+2\cdot 4(3-\frac{\nu}{n})n^2=2(\nu-2n)^2+16n^2.
$$
Therefore, since $\mathop{\rm d}(Z)=4n^2$, the inequality
$$
\frac{\mathop{\rm mult}\nolimits_o}{\mathop{\rm deg}} Z>\frac{4}{M}
$$
holds, which contradicts the estimate (\ref{04.11.2017.1}) for
$a=2$, $c=1$ (see Remark 2.1 (i)). This contradiction proves Proposition 3.3. Q.E.D.\vspace{0.1cm}

So the codimension $\alpha$ can take at most two values: 2 and 3. In order to exclude these options as well, let us consider one more characteristic of the subvariety $B$. For a pair of distinct points
$p,q\in B$ the symbol $[p,q]$ denotes the line in $E_{\mathbb P}$, joining the points $p$ and $q$, {\it provided that this line is contained in} $Q$. If this line is not contained in $Q$, then set $[p,q]=\emptyset$. Now set:
$$
\mathop{\rm Sec} (B\subset Q)=\overline{\mathop{\bigcup}\limits_{(p,q)\in
B\times B\setminus\Delta_B} [p,q]},
$$
where $\Delta_B\subset B\times B$ is the diagonal. Since
$B\subset \mathop{\rm Sec} (B\subset Q)$, we have
$$
\mathop{\rm codim}(\mathop{\rm Sec}(B\subset Q)\subset Q)\leqslant\mathop{\rm codim} (B\subset Q).
$$
We say that the case $\alpha.\beta$ takes place, where
$\alpha\in\{2,3\}$ and $0\leqslant \beta\leqslant\alpha$, if
$\mathop{\rm codim}(B\subset Q)=\alpha$ and
$$
\mathop{\rm codim}(\mathop{\rm Sec}(B\subset Q)\subset Q)=\beta.
$$

{\bf Proposition 3.4.} {\it Let $X\subset{\mathbb P}^N$ be an irreducible subvariety of codimension 2 and $N\geqslant 3$. Then
$\mathop{\rm Sec}(X)=\langle X\rangle$, that is, one of the three options takes place:

\begin{itemize}

\item $\mathop{\rm Sec}(X)={\mathbb P}^N$,

\item $X$ is a hypersurface of degree $d_X\geqslant 2$ in some hyperplane in ${\mathbb P}^N$,

\item $X$ is a linear subspace of codimension two.

\end{itemize}
}

{\bf Proof:} this is obvious. Q.E.D.\vspace{0.1cm}

In the case $\alpha=2$ we apply Proposition 3.4 to the intersection
$B\cap\Pi$, where $\Pi\subset Q$ is a general linear subspace of maximal dimension. We get the following list of options. For a proof, see Subsection 4.2.\vspace{0.1cm}

{\bf Case 2.0.} $\mathop{\rm Sec} (B\subset Q)=Q$.\vspace{0.1cm}

{\bf Case 2.1.} $\mathop{\rm Sec} (B\subset Q)$ is a hyperplane section of the quadric $Q$, on which the subvariety $B$ is cut out by a hypersurface of degree $\mathop{\rm Sec} (B)\geqslant 2$.\vspace{0.1cm}

{\bf Case 2.2.} $\mathop{\rm Sec} (B\subset Q)=B$ is a section of the quadric $Q$ by a linear subspace of codimension 2.\vspace{0.1cm}

For $\alpha=3$ we can claim the following.\vspace{0.1cm}

{\bf Case 3.0.} $\mathop{\rm Sec} (B\subset Q)=Q$.\vspace{0.1cm}

{\bf Case 3.1.} $\mathop{\rm Sec} (B\subset Q)$ is a divisor on the quadric
$Q$.\vspace{0.1cm}

{\bf Case 3.2.} $\mathop{\rm Sec} (B\subset Q)$ is a section of the quadric
$Q$ by a linear subspace of codimension 2, on which the subvariety
$B$ is cut out by a hypersurface of degree $\mathop{\rm Sec} (B)\geqslant 2$.\vspace{0.1cm}

{\bf Case 3.3.} $\mathop{\rm Sec} (B\subset Q)=B$ is a section of the quadric $Q$ by a linear subspace of codimension 3.\vspace{0.1cm}

In the description of the cases $3.\beta$ given above, the only not quite obvious statement is the description of the case 3.2, which is based on an analog of Proposition 3.4 for a subvariety  of codimension 3. This analog is stated and proven in Section 4 (the proof is very simple).\vspace{0.1cm}

{\bf Remark 3.1.} Let us sum up what has been done so far. Assuming the existence of a maximal singularity $E^*$ of the linear system
$\Sigma$ (and the non-existence of maximal subvarieties of the form
$P\cap V$, where $P\subset{\mathbb P}$ is a linear subspace of codimension 2), we proved the existence of another singularity
$E$ of the linear system $\Sigma$, the centre of which is a double point
$o\in V$. This singularity $E$ ``looks like a maximal one'' in the sense that it satisfies the {\it Noether-Fano type inequality}
\begin{equation}\label{15.11.2017.1}
\mathop{\rm ord}\nolimits_E\Sigma>
n\cdot (3\mathop{\rm ord}\nolimits_E Q+a(E,V^+))
\end{equation}
(in the brackets we could have added $+1$, since the pair $\Box$ is not log canonical, but we do not need this). That inequality is weaker than the standard Noether-Fano inequality, but this is compensated by a high dimension of the centre $B$ of the singularity $E$ on $V^+$. We will show below that in each of the cases $\alpha.\beta$ listed above, the existence of the singularity $E$ leads to a contradiction. Thus, the existence of the original maximal singularity $E^*$ leads to a contradiction, either. This will complete the proof of Theorem 5.\vspace{0.3cm}


{\bf 3.3. Resolution of the singularity $E$.} Let us consider the standard procedure of resolving the singularity $E$. Let $\varphi_{i,i-1}\colon V_i\to V_{i-1}$ be the blow up of the centre $B_{i-1}$ of the singularity $E$ on $V_{i-1}$, where $V_0=V$, so that $B_0=o$ and $V_1=V^+$. The exceptional divisor of the blow up $\varphi_{i,i-1}$ is denoted by the symbol $E_i$, so that $E_1=Q$. The sequence of these blow ups terminates: for some $i=K$ the exceptional divisor $E_K$ is the centre of $E$ on $V_K$, and there is nothing more to blow up. The varieties $V_i$ have, generally speaking, uncontrollable singularities, however, at the general point of the subvariety $B_i$ the variety $V_i$ is non-singular for $i\geqslant 1$, and this is the only property that we need for all computations. For $i>j$ the composition of blow ups
$$
V_i\to V_{i-1}\to\cdots \to V_j
$$
is denoted by the symbol $\varphi_{i,j}$. For all details, see
\cite[Chapter 2, Subsection 1.2]{Pukh13a}.\vspace{0.1cm}

On the set $1,\dots, K$ of indices, numbering the exceptional divisors, we introduce the structure of an oriented graph, setting $i\to j$, if
$i>j$ and the centre of the blow up $B_{i-1}$ is contained in the strict transform $E_j^{i-1}\subset V_{i-1}$ of the exceptional divisor $E_j$ on
$V_{i-1}$. In particular, we always have $(i+1)\to i$. By the symbol $p_{ij}$ we denote the number of paths in the just constructed oriented graph from the vertex $i$ to the vertex $j$, if $i\neq j$; by definition, $p_{ii}=1$. For $i<j$, obviously, $p_{ij}=0$. In order to simplify the notations, we write $p_i$ instead of $p_{Ki}$. For $i=1,\dots, K$ we define the numbers $\nu_i\in{\mathbb Z}_+$ (the ``elementary multiplicities''), taking a general divisor $D\in\Sigma$ and writing down
$$
D^i=\varphi^*_{i,i-1}(D^{i-1})-\nu_iE_i,
$$
where the upper index $j$ denotes, as usual, the operation of taking the strict transform on $V_j$ (see \cite{Pukh13a}). Now
$$
\mathop{\rm ord}\nolimits_E\Sigma=\sum^K_{i=1}p_i\nu_i
$$
and
$$
\mathop{\rm ord}\nolimits_E Q=p_1,\quad a(E,V^+)=\sum^K_{i=2}p_i\delta_i,
$$
where $\delta_i=\mathop{\rm codim}(B_{i-1}\subset V_{i-1})-1$ is the
``elementary discrepancy''. Now the Noether-Fano type inequality
(\ref{15.11.2017.1}) can be explicitly re-written in the form of the estimate
\begin{equation}\label{15.11.2017.2}
\sum^K_{i=1}p_i\nu_i>n\cdot\left(3p_1+\sum^K_{i=2}p_i\delta_i\right).
\end{equation}
Since from the definition of the numbers $p_i=p_{Ki}$ we have the obvious equality $a>i$:
$$
p_{ai}=\sum_{a\to j}p_{ji},
$$
we conclude that if $\nu_K\leqslant \delta_Kn$, then for some
$K_1<K$ the estimate
$$
\sum^{K_1}_{i=1}p_{K_1i}\nu_i>n\cdot\left(3p_{K_11}+
\sum^{K_1}_{i=2}p_{K_1i}\delta_i\right)
$$
holds. For this reason we may (and will) assume that $\nu_K>\delta_Kn$. Furthermore, we break the set of vertices $\{1,\dots, K\}$ of the constructed oriented graph into the {\it lower part} $\{1,\dots, L\}$ and the {\it upper part} (which can be empty) $\{L+1,\dots, K\}$, setting $i\leqslant L$, if and only if $\delta_i\geqslant 2$. Finally, we use the well known procedure of
{\it erasing arrows} in the constructed graph: let us remove all arrows
$i\to 1$, going from the vertices of the upper part ($i\geqslant L+1$). Recall that by Proposition 3.3 we have $L\geqslant 2$, so that at least the vertex 2 lies in the lower part. The procedure of erasing arrows may decrease $p_1$, but does not change the numbers $p_i$ for $i\geqslant 2$, therefore the inequality (\ref{15.11.2017.2}) can only get stronger. For this reason we assume that in the oriented graph there are no arrows from the vertices of the upper part to the vertex 1. Set
$$
\Sigma^+_0=\sum_{i\geqslant 2, \delta_i=3}p_i,\quad
\Sigma_1=\sum_{\delta_i=2}p_i,\quad \Sigma_2=\sum_{\delta_i=1}p_i
$$
and $\Sigma_0=p_1+\Sigma^+_0$. The procedure of erasing arrows gives the following fact.\vspace{0.1cm}

{\bf Proposition 3.5.} {\it The following inequality holds:}
$$
p_1\leqslant \Sigma^+_0+\Sigma_1.
$$

{\bf Proof.} Indeed, every path from the vertex $K$ to the vertex 1 is of the form
$$
K\to\cdots\to a\to 1,
$$
where $a\leqslant L$. Therefore,
$$
p_1=p_{K1}=\sum_{a\to 1}p_{Ka}\leqslant \Sigma^+_0+\Sigma_1,
$$
which what we need. Q.E.D. for the proposition.\vspace{0.1cm}

Now let us consider the self-intersection $Z=(D_1\circ D_2)$ of the system
$\Sigma$. Set
$m_1=\mathop{\rm m}(Z)=\frac12 \mathop{\rm mult}\nolimits_o Z$
and for $i=2,\dots, L$
$$
m_i=\mathop{\rm mult}\nolimits_{B_{i-1}}Z^{i-1}
$$
(the definition makes sense, because for $i\leqslant L$ the subvariety
$B_{i-1}\subset V_{i-1}$ is of codimension
$\geqslant 3$, whereas the strict transform $Z^{i-1}$ of the cycle
$Z$ on $V_{i-1}$ is an effective cycle of codimension 2). Writing down
$$
(D^1_1\circ D^1_2)=Z^1+Z_1,
$$
where $Z_1$ is an effective divisor on the quadric $Q=E_1$, we obtain the equality
$$
m_1=\nu^2+\mathop{\rm d}\nolimits_Q(Z_1),
$$
whence, arguing as in \cite{Pukh2017a} and in \cite[Chapter 2]{Pukh13a}, we obtain by means of the technique of counting multiplicities the following estimate:
\begin{equation}\label{15.11.2017.3}
\sum^L_{i=1}p_im_i\geqslant \sum^K_{i=1}p_i\nu^2_i,
\end{equation}
where $m_1\geqslant m_2\geqslant\cdots \geqslant m_L$ and
$\nu_1\geqslant \nu_2\geqslant\cdots \geqslant \nu_K$. (Note that the first inequalities in both expressions, that is, $m_1\geqslant m_2$ and
$\nu_1\geqslant\nu_2$, are non-trivial, although for the quadric $Q$ their proof is very simple; for the general case see \cite[Proposition 2]{Pukh2017a}). If there is no additional information about the multiplicities $\nu_i$, then we estimate the minimum of the quadratic form in the right hand side of the inequality (\ref{15.11.2017.3}) on the hyperplane, the equation of which is obtained from the inequality (\ref{15.11.2017.2}) by replacing the $>$ sign by $=$, and get the estimate
$$
\sum^L_{i=1}p_im_i>\frac{(3p_1+3\Sigma_0^++2\Sigma_1+\Sigma_2)^2}{p_1+
\Sigma_0^++\Sigma_1+\Sigma_2}n^2.
$$
Replacing the left hand side of this inequality by
$$
p_1m_1+m_2\sum^L_{i=1}p_i
$$
(which can only make the inequality sharper), we get finally:
\begin{equation}\label{15.11.2017.4}
p_1m_1+(\Sigma^+_0+\Sigma_1)m_2>\frac{(3p_1+3\Sigma_0^++2\Sigma_1+
\Sigma_2)^2}{p_1+\Sigma_0^++\Sigma_1+\Sigma_2}n^2.
\end{equation}
Estimating the left hand side from above by the expression $(\Sigma_0+\Sigma_1)p_1$, we get the $4n^2$-inequality $m_1>4n^2$, mentioned above (as it was done in \cite{Pukh2017a}). However, we can say more.\vspace{0.1cm}

{\bf Proposition 3.6.} {\it The following inequality holds:}
$$
m_1+m_2>8n^2.
$$

{\bf Proof.} For a fixed value of the linear form
$p_1m_1+(\Sigma^+_0+\Sigma_1)m_2$ (in the variables
$m_1\geqslant m_2\geqslant 0$) the minimum of the expression $m_1+m_2$ is attained for $m_1=m_2$ (recall that $p_1\leqslant\Sigma^+_0+\Sigma_1$). In that case $m_1=m_2>4n^2$, which proves Proposition 3.6. Q.E.D.\vspace{0.1cm}

For certain types of the singularity $E$ there is some additional information about the multiplicities $m_i$ and $\nu_i$, which makes it possible to make the estimates for $m_1$ and $m_2$ sharper. In order to obtain such information, we need some facts about the secant variety $\mathop{\rm Sec}(B\subset Q)$. Thses questions are dealt with in the next section.


\section{Subvarieties of the quadric $Q$}

The aim of this section is to prove the classification of options given in Subsection 3.2. First we consider the problem of irreducibility of the intersection of an irreducible subvariety $X\subset Q$ of a small codimension with a general linear subspace $P\subset Q$ of maximal dimension (Subsection 4.1). On this basis it is easy to prove our classification, that is, the description of the cases 2.1, 2.2 and 3.2, 3.3 (Subsection 4.2). In Subsection 4.3 we discuss non-degenerate subvarieties $X\subset{\mathbb P}^N$ of codimension 3, the secant variety $\mathop{\rm Sec}  (X)$ of which are strictly smaller than ${\mathbb P}^N$.\vspace{0.3cm}

{\bf 4.1. Irreducibility of the intersection with a linear subspace.} Let
${\cal L}$ be a closed algebraic set, parameterizing linear subspaces of maximal dimension on the quadric $Q$, that is, of dimension
$M-\ulcorner\frac12 \mathop{\rm rk}Q\urcorner$. Depending on whether the rank of the quadric is odd or even, the set ${\cal L}$ can be irreducible variety or a union of two irreducible varieties. For a subspace $P\in{\cal L}$ by the symbol $\pi_P$ we denote the projection ${\mathbb P}^M$ from the subspace $P$. For a subspace $\Lambda\subset{\mathbb P}^M$, such that $\Lambda\supset P$ and $\mathop{\rm dim} \Lambda= \mathop{\rm dim} P+1$ (that is, a fibre of the projection $\pi_P$) we have $Q\cap\Lambda=P\cup Q(P,\Lambda)$, where
$Q(P,\Lambda)\in{\cal L}$. If $P\in{\cal L}$ is a subspace of general position and $\Lambda$ a general fibre of the projection $\pi_P$, then
$Q(P,\Lambda)$ is also a subspace of general position.\vspace{0.1cm}

{\bf Proposition 4.1.} {\it Let $X\subset Q$ be an irreducible subvariety of codimension $\leqslant 3$ and $P\in{\cal L}$ a linear subspace of general position. Then $X\cap P$ is an irreducible subvariety of codimension
$\mathop{\rm codim}  (X\subset Q)$ in the projective space
$P$.}\vspace{0.1cm}

{\bf Proof.} The equality
$$
\mathop{\rm codim}  ((X\cap P)\subset P)=\mathop{\rm codim}  (X\subset Q)
$$
is obvious, we only need to show the irreducibility. Assume that
$$
X\cap P=\mathop{\bigcup}\limits_{i\in I}X_i(P),
$$
where $|I|\geqslant 2$. Then for any general $P\in{\cal L}$ we have a similar picture with the same value of $|I|$. For a fibre $\Lambda$ of the projection $\pi_P$ the intersection
$$
P_{\Lambda}=P\cap Q(P,\Lambda)
$$
is a hyperplane in $P$, and moreover, it is easy to check that varying $\Lambda$, we obtain the complete family of hyperplanes in $P$, containing the vertex space of the quadric $Q$. For that reason for a general space
$\Lambda\supset P$ the closed space
$$
X_i(P)\cap P_{\Lambda}
$$
is irreducible for all $i\in I$, and we may assume that all sets
$X_i(P)\cap P_{\Lambda}$, $i\in I$, are distinct, so that the components
$X_i(P)$ are identified by the intersections with the hyperplane
$P_{\Lambda}$. However, our construction is symmetric with respect to the subspaces $P$ and $Q(P,\Lambda)$, so that for a general fibre
$\Lambda$ of the projection $\pi_P$ there is a bijective correspondence between irreducible components of the intersection $X\cap P$ and the components of the intersection $X\cap Q(P,\Lambda)$, which makes it possible to write down
$$
X\cap Q(P,\Lambda)=\mathop{\bigcup}\limits_{i\in I}X_i(Q(P,\Lambda)).
$$
But then for each $i\in I$
$$
X_i=\overline{\mathop{\bigcup}\limits_{\Lambda\in U} X_i(Q(P,\Lambda))}
$$
(the union is taken over a non-empty Zariski open subset $U\subset {\mathbb P}^{M-1-\mathop{\rm dim} P}$) is an irreducible component of the original set $X$, which contradicts the assumption about its irreducibility. Q.E.D. for the proposition.\vspace{0.1cm}

{\bf Remark 4.4.} We stated and proved Proposition 4.1 in the form in which it will be used. The proof given above works word for word for the case of an irreducible subvariety of codimension
$[\frac12 \mathop{\rm rk} Q-1]$.\vspace{0.3cm}


{\bf 4.2. Secant subvarieties on the quadric $Q$.} We will not consider the problem of reducibility of the secant subvariety $\mathop{\rm Sec}  (B\subset Q)$. In order to prove the classification of options given in Subsection 3.2, it is sufficient to consider the component ${\mathop{\rm Sec} }^* (B\subset Q)$, defined by the following condition: for a general subspace
$P\in {\cal L}$ the set ${\mathop{\rm Sec} }^* (B\subset Q)$ contains the secant variety
$$
\mathop{\rm Sec}  (B\cap P)\subset P.
$$
Obviously, $\mathop{\rm dim}{\mathop{\rm Sec} }^* (B\subset Q)=\mathop{\rm dim}\mathop{\rm Sec}  (B\subset Q)$. In order not to make the notations too complicated, by the secant variety $\mathop{\rm Sec}  (B\subset Q)$ we will mean that very component of the maximal dimension.\vspace{0.1cm}

Now let us consider the classification of possible cases given in Subsection 3.2. Assume first that $\mathop{\rm codim}  (B\subset Q)=2$. Then for a general subspace $P\in{\cal L}$ we have
$$
\mathop{\rm codim} ((B\cap P)\subset P)=2,
$$
where $B\cap P$ is an irreducible subvariety of degree
$\mathop{\rm d}_Q (B)$ in the projective space $P$. If
$\mathop{\rm Sec} (B\subset Q)$ has codimension 2, then by Proposition 3.4 $B\cap P$ is a linear subspace of codimension 2 in
$P$. Therefore, $\mathop{\rm d}_Q (B)=1$ and $\mathop{\rm deg} B=2$. But then $B$ is contained in a linear subspace
$\Pi\subset{\mathbb P}^M$ of dimension $\mathop{\rm dim} B+1=M-2$ and for that reason $B=Q\cap \Pi$, just as the case 2.2 claims.\vspace{0.1cm}

If $\mathop{\rm Sec} (B\subset Q)$ is of codimension 1 in $Q$, then for a general subspace $P\in{\cal L}$
$$
\mathop{\rm codim} ((B\cap P)\subset P)\geqslant 1.
$$
Taking into account the previous case, we can claim that the last inequality is an equality. According to Proposition 3.4, in this case
$$
\mathop{\rm Sec} (B\subset Q)\cap P
$$
contains a hyperplane in $P$, and for that reason by Proposition 4.1 is a hyperplane in $P$. Then $\mathop{\rm d}_Q (\mathop{\rm Sec} (B\subset Q))=1$, so that $\mathop{\rm deg} \mathop{\rm Sec} (B\subset Q)=2$. Therefore, $\mathop{\rm Sec} (B\subset Q)$ is a hyperplane section of the quadric $Q$, which is itself a factorial quadric. Since
$\mathop{\rm d}_Q (B)\geqslant 2$ (if $\mathop{\rm d}_Q (B)=1$, then we are in the case 2.2), we obtain precisely the description of the case 2.1. The description of the case 2.0 requires no proof.\vspace{0.1cm}

Assume now that $\mathop{\rm codim} (B\subset Q)=3$. If
$\mathop{\rm Sec} (B\subset Q)$ has codimension 3, then we argue as above in the case 2.2 and obtain the description of the case 3.3. If
$\mathop{\rm Sec} (B\subset Q)$ has codimension 2, then we argue as above in the case 2.1, using the following simple fact.\vspace{0.1cm}

{\bf Proposition 4.2.} {\it Let $X\subset{\mathbb P}^N$ be an irreducible subvariety of dimension $\mathop{\rm dim} X\geqslant 2$, and moreover,
$\mathop{\rm dim} \mathop{\rm Sec} (X)=\mathop{\rm dim} X+1$. Then
$\mathop{\rm Sec} (X)$ is a linear subspace and $X$ is a hypersurface of degree $d_X\geqslant 2$ in that subspace.}\vspace{0.1cm}

{\bf Proof.} For a point $p\in X$ of general position consider the cone
$C(p,X)$ with the vertex $p$ and the base $X$ (the closure of the union of all lines $[p,q]$, where $q\in X\setminus \{p\}$). This is an irreducible subvariety of some degree $d\geqslant 1$ and dimension
$\mathop{\rm dim} X+1$, so that $\mathop{\rm Sec} (X)=C(p,X)$. Therefore,
$$
\mathop{\rm mult}\nolimits_p \mathop{\rm Sec} (X)=d=\mathop{\rm deg}\mathop{\rm Sec}  (X)
$$
for all points $p\in X$, that is, for every point $p\in X$ the subvariety
$\mathop{\rm Sec} (X)$ is a cone with the vertex at the point $p$. Since
$X$ is a subvariety of codimension 1 in $\mathop{\rm Sec} (X)$, it is easy to see that $d=1$. Q.E.D. for the proposition.\vspace{0.1cm}

The proposition proven above implies the description of the case 3.2. The cases 3.1 and 3.0 require no proof: in those two cases we just note the codimension of the set $\mathop{\rm Sec} (B\subset Q)$. This completes the proof of the list of options given in Subsection 3.2.\vspace{0.3cm}


{\bf 4.3. A remark on the secant variety.} In
\cite[Section 3]{Pukh16a} it was shown that if
$B\subset V$ is a subvariety of codimension 2 on a general smooth hypersurface $V\subset{\mathbb P}$ of degree $M$ and $B$ is not contained in a hyperplane (that is, $\langle B\rangle={\mathbb P}$), then
$\mathop{\rm Sec} (B)={\mathbb P}$. Since in that case
$\mathop{\rm codim}  (B\subset{\mathbb P})=3$, it is natural to ask: what is an analog of Proposition 3.4 for subvarieties of codimension 3?\vspace{0.1cm}

Let $X\subset{\mathbb P}^N$ be an irreducible subvariety of codimension 3. If $\mathop{\rm dim} \langle X\rangle\leqslant N-1$, then applying Proposition 3.4 to the projective space $\langle X\rangle$, we get a complete classification of options. Assume therefore that
$\langle X\rangle={\mathbb P}^N$.\vspace{0.1cm}

{\bf Example 4.1.} Let $\Gamma\subset{\mathbb P}^4$ be a non-degenerate curve. Obviously, $\mathop{\rm Sec} (\Gamma)$ is a hypersurface in
${\mathbb P}^4$. Considering the cone over
$\Gamma\subset{\mathbb P}^4\subset{\mathbb P}^N$ with a vertex subspace of dimension $N-5$, we obtain an irreducible subvariety
$X\subset{\mathbb P}^N$, for which $\mathop{\rm Sec} (X)$ is a hypersurface in ${\mathbb P}^N$. We say that a subvariety $X$, obtained in this way, is a {\it cone over a non-degenerate curve in} ${\mathbb P}^4$. The following fact takes place.\vspace{0.1cm}

{\bf Proposition 4.3.} {\it Let $X\subset {\mathbb P}^N$ be a non-degenerate irreducible subvariety of codimension 3 and assume that
$$
\mathop{\rm Sec} (X)\neq{\mathbb P}^N.
$$
Then $X$ is a cone over a non-degenerate curve in} ${\mathbb P}^4$.\vspace{0.1cm}

We do not give a {\bf proof} here, because we do not use this fact. The proof follows the same scheme of arguments as in \cite[Section 3]{Pukh16a}. A close look at that proof shows that the condition $B\subset V$ was used only to claim that for every point $p\in B$ there are at most finitely many lines on $B$ passing through this point.\vspace{0.1cm}

The proof given in \cite[Section 3]{Pukh16a} can be improved to a proof of Proposition 4.3.


\section{Exclusion of the maximal singularity}

In this section we prove Theorem 5. First we exclude the case 3.3 (Subsection 5.1), then the case 3.2 (Subsection 5.2). These are the most difficult cases, requiring considerable efforts; the remaining five cases are excluded in Subsection 5.3. This completes the proof of Theorem 5.\vspace{0.3cm}

{\bf 5.1. Exclusion of the case 3.3.} First of all, let us recall certain notations which we will use in this section (in Subsections 5.1 -- 5.3). There is a fixed point $o\in V$, a quadratic singularity of rank
$\geqslant 13$. Its blow up in the projective space ${\mathbb P}$ is denoted by the symbol
$\varphi_{{\mathbb P}}\colon {\mathbb P}^+\to{\mathbb P}$, and the exceptional divisor $\varphi^{-1}_{{\mathbb P}}(o)$ --- by the symbol $E_{\mathbb P}$. The strict transform on ${\mathbb P}^+$ is denoted by adding the upper index $+$, for instance, $\varphi\colon V^+\to V$ is the blow up of the point $o$ on the variety $V$, where $Q=E_{\mathbb P}\cap V^+$ is the exceptional quadric. Furthermore, $B\subset Q$ is the centre of the singularity $E$ on $V^+$ (see \S 3). The self-intersection $(D_1\circ D_2)$ of the mobile linear system $\Sigma\subset |2nH|$ is denoted by the symbol $Z$. Hypertangent divisor at the point $o$ are denoted by the symbols $T_i$, where $T_i\in\Lambda_i$ (see \S 2).\vspace{0.1cm}

Assume that the case 3.3 takes place. Our aim is to obtain a contradiction. By Proposition 3.6 the inequality
$$
\mathop{\rm m}(Z)+\mathop{\rm mult}\nolimits_B Z^+>2\mathop{\rm d} (Z)
$$
holds. This inequality is linear in $Z$. For that reason we may assume that $Z$ is an irreducible subvariety of codimension 2. By Proposition 2.4 we have $\mathop{\rm m}(Z)\leqslant 2\mathop{\rm d} (Z)$, so that
$Z^+$ contains $B$ and thus $Z\not\subset T_2$ (by the regularity condition (R2.1)).\vspace{0.1cm}

Set $\overline{\Pi}\subset {\mathbb P}$ to be the uniquely determined subspace of codimension 3, containing the point $o$ and ``cutting out
$B$ on $Q$'', that is,
$$
\overline{\Pi}^+\cap Q=B.
$$
The corresponding section $V\cap\overline{\Pi}$ is denoted by the symbol
$\Delta$ (in Subsection 2.4 it was denoted by the symbol $V_{\Pi}$). So
$\Delta^+\cap Q=B$. Consider the linear system $|H-\Delta|$ of hyperplane sections of the hypersurface $V$, containing $\Delta$. Let
$R_1\in|H-\Delta|$ be a general divisor. Since
$B\not\subset T_2^+$, we have $\Delta\not\subset T_2$, which implies that none of the irreducible components of the effective cycle $(Z\circ R_1)$ of codimension 3 is contained in $T_2$.\vspace{0.1cm}

Since $\mathop{\rm Bs} |H-\Delta|=\Delta$, and
$\mathop{\rm Bs} |H-\Delta|^+=\Delta^+$ (in the scheme-theoretic sense), the following equalities hold:
$$
\mathop{\rm m} (Z\circ R_1)=\mathop{\rm m} (Z),\quad \mathop{\rm mult}\nolimits_B(Z\circ R_1)^+=\mathop{\rm mult}\nolimits_B Z^+.
$$
Besides,
$\mathop{\rm d}(Z\circ R_1)=\mathop{\rm d} (Z)$. Write down
$$
(Z\circ R_1)=a\Delta+Z_{\sharp},
$$
where $a\in {\mathbb Z}_+$ and the effective cycle $Z_{\sharp}$ does not contain $\Delta$ as a component. Since
$\mathop{\rm d}(\Delta)=1$, $\mathop{\rm m}(\Delta)=1$ and
$\mathop{\rm mult}\nolimits_B\Delta^+=1$, we obtain the inequality
$$
\mathop{\rm m}(Z_{\sharp})+\mathop{\rm mult}\nolimits_B Z_{\sharp}^+ >2\mathop{\rm d}(Z_{\sharp}).
$$
Let us consider one more general divisor $R_2\in |H-\Delta|$. Obviously, none of the components of the cycle $Z_{\sharp}$ is contained in $R_2$, so that the effective cycle $(Z_{\sharp}\circ R_2)$ is well defined. It has codimension 4 in $V$, 3 in $R_1$ and 2 in $R_1\cap R_2$. The following inequality holds:
$$
\mathop{\rm m}(Z_{\sharp}\circ R_2)\geqslant
\mathop{\rm m}(Z_{\sharp})+\mathop{\rm mult}\nolimits_B Z_{\sharp}^+>2\mathop{\rm d}(Z_{\sharp}).
$$
Since
$\mathop{\rm d}(Z_{\sharp})=\mathop{\rm d}(Z_{\sharp}\circ R_2)$,
there is an irreducible component $Y$ of the cycle $(Z_{\sharp}\circ R_2)$, satisfying the inequality
\begin{equation}\label{22.11.2017.1}
\mathop{\rm m}(Y)>2\mathop{\rm d}(Y).
\end{equation}

{\bf Lemma 5.1.} {\it $Y$ is not contained in} $T_2$.\vspace{0.1cm}

{\bf Proof.} Assume the converse: $Y\subset T_2$. By construction,
$Y$ is an irreducible component of the intersection of the divisor $R_2$ with one of the irreducible components of the cycle $Z_{\sharp}$, which, as we know, is not contained in $T_2$. Therefore, $Y$ is an irreducible component of the effective cycle $(Z_{\sharp}\circ T_2)$, {\it which is contained in} $R_2$. The cycle $(Z_{\sharp}\circ T_2)$ is of codimension 3 in $R_1$. Since that cycle does not depend on $R_2$, and $R_2\in |H-\Delta|$ by assumption is a general divisor of this linear system, we conclude that $Y\subset \Delta$.\vspace{0.1cm}

So $Y$ is a prime divisor on $\Delta$. By the condition
(R2.1) the divisor $(T_2\circ\Delta)$ on $\Delta$ is irreducible and reduced, so that we obtain the equality $Y=(T_2\circ\Delta)$. But this is impossible: by the condition (R2.1) we have $\mathop{\rm m}(T_2\circ\Delta)=3$. Since
$\mathop{\rm d}(T_2\circ\Delta)=2$, we get a contradiction with the inequality (\ref{22.11.2017.1}). Q.E.D. for the lemma.\vspace{0.1cm}

Now let us consider the effective cycle $(Y\circ T_2)$ of codimension 3 on
$R_1\cap R_2)$. It satisfies the inequality
$$
\frac{\mathop{\rm mult}\nolimits_o}{\mathop{\rm deg}} (Y\circ T_2)>\frac{6}{M},
$$
which contradicts the inequality (\ref{04.11.2017.1}) for $c=2$ (Proposition 2.4), since
$$
\frac{6}{M}\geqslant \frac{5}{M-1}
$$
already for $M\geqslant 6$. This contradiction completes the exclusion of the case 3.3.\vspace{0.3cm}


{\bf 5.2. Exclusion of the case 3.2.} Assume that the case 3.2 takes place. Again our aim is to obtain a contradiction. Now by the symbol $\Delta$ we denote the section of the hypersurface $V$ by the uniquely determined subspace $\overline{\Pi}\subset{\mathbb P}$ of codimension 2, such that
$$
\overline{\Pi}^+ \cap Q=\mathop{\rm Sec} (B\subset Q).
$$
(Recall that $\mathop{\rm Sec} (B\subset Q)$ is the section of the quadric
$Q\subset E_{\mathbb P}={\mathbb P}^M$ by the linear subspace
$\langle B\rangle\subset{\mathbb P}^M$ of codimension 2.) Therefore,
$$
\Delta^+\cap Q=\mathop{\rm Sec} (B\subset Q).
$$
We continue to assume that the self-intersection $Z$ is an irreducible subvariety of codimension 2. Since $\mathop{\rm d}(\Delta)=1$ and
$\mathop{\rm m}(\Delta)=\mathop{\rm mult}\nolimits_B\Delta^+=1$, we have
$Z\neq\Delta$. Consider the pencil of hyperplane sections $|H-\Delta|$, containing $\Delta$, and take a general divisor $R\in|H-\Delta|$. By construction, $Z\not\subset R$, so that the effective cycle
$Z_R=(Z\circ R)$ of codimension 2 on the hyperplane section $R$ is well defined. Setting for the convenience of notations $\Delta_Q=\Delta^+\cap Q$, let us write down
$$
(Z^+\circ R^+)=Z^+_R+a\Delta_Q,
$$
where $a\in{\mathbb Z}_+$. Such a writing is possible, because
$\Delta_Q=\mathop{\rm Bs} |H-\Delta|^+_Q$. Now we get
$\mathop{\rm m}(Z_R)=\mathop{\rm m}(Z)+a$ and
$a+\mathop{\rm mult}\nolimits_B Z^+_R\geqslant \mathop{\rm mult}\nolimits_B Z^+$, since obviously $\mathop{\rm mult}\nolimits_B \Delta_Q=1$. From there we obtain the inequality
\begin{equation}\label{23.11.2017.1}
\mathop{\rm m}(Z_R)+\mathop{\rm mult}\nolimits_B Z^+_R>2\mathop{\rm d}(Z)=2\mathop{\rm d}(Z_R).
\end{equation}
By the linearity of this inequality in $Z_R$ we may assume that
$Z_R=Y$ is an irreducible subvariety of codimension 2 (with respect to
$R$).\vspace{0.1cm}

We know that $B$ is a prime divisor on the quadric
$\Delta_Q\subset\langle B\rangle={\mathbb P}^{M-2}$, cut out on this quadric by a hypersurface of degree $\mathop{\rm d}_Q(B)\geqslant 2$. The effective cycle $(Y^+\circ Q)$ of dimension $\mathop{\rm dim} B$ contains
$B$ with multiplicity $\geqslant \mathop{\rm mult}\nolimits_B Y^+$. Since
$$
\mathop{\rm m}(Y)=\mathop{\rm d}\nolimits_Q(Y^+\circ Q),
$$
this implies the inequality
$$
\mathop{\rm m}(Y)\geqslant \mathop{\rm d}\nolimits_Q(B)\cdot \mathop{\rm mult}\nolimits_B Y^+.
$$
Now we have to consider two cases:\vspace{0.1cm}

(1) $Y\subset\Delta$ is a prime divisor,\vspace{0.1cm}

(2) $Y$ is not contained in $\Delta$.\vspace{0.1cm}

{\bf Lemma 5.2.} {\it In the case (1) the equality
$\mathop{\rm d}_Q(B)=2$ holds.}\vspace{0.1cm}

{\bf Proof.} Assume the converse:
$\mathop{\rm d}_Q(B)\geqslant 3$. We have the inequality
$\mathop{\rm m}(Y)\geqslant 3\mathop{\rm mult}\nolimits_B Y^+$. Combining it with the inequality (\ref{23.11.2017.1}) for $Z_R=Y$, we obtain the estimate
$$
\mathop{\rm m}(Y)>\frac32 \mathop{\rm d}(Y).
$$
From the regularity conditions we deduce that
$Y\neq (T_2\circ\Delta)$. Indeed, the quadratic form
$q_2|_{\Pi}$ is of rank $\geqslant 9$, the variety $\Delta$ is factorial, so that the divisor $(T_2\circ\Delta)$ is irreducible and reduced. By the condition (R2.1) this divisor has multiplicity precisely 6 at the point $o$ and for that reason the equality
$$
\mathop{\rm m}(T_2\circ\Delta)=\frac32 \mathop{\rm d}(T_2\circ\Delta)
$$
holds, from which we get that $Y\neq (T_2\circ\Delta)$). So we obtain the following well defined effective cycle
$$
(Y\circ(T_2\circ\Delta))
$$
of codimension 2 on $\Delta$, satisfying the inequality
$$
\mathop{\rm m}(Y\circ (T_2\circ\Delta))>\frac94
\mathop{\rm d}(Y\circ (T_2\circ\Delta)),
$$
which can be re-written in the form
$$
\frac{\mathop{\rm mult}\nolimits_o}{\mathop{\rm deg}} (Y\circ (T_2\circ\Delta))>
\frac{9}{2M}.
$$
This contradicts the inequality (\ref{04.11.2017.1}) for
$a=c=2$, since for $M\geqslant 9$
$$
\frac{9}{2M}\geqslant \frac{4}{M-1}.
$$
The just obtained contradiction proves Lemma 5.2.\vspace{0.1cm}

{\bf Lemma 5.3.} {\it The case (2) is impossible.}\vspace{0.1cm}

{\bf Proof.} Assume the converse: the case (2) takes place. Since $\mathop{\rm d}_Q(B)\geqslant 2$, the estimate
$\mathop{\rm m}(Y)\geqslant 2\mathop{\rm mult}\nolimits_B Y^+$ holds. Let us consider the effective cycle $(Y\circ\Delta)=(Y\circ\Delta)_R$ of the scheme-theoretic intersection of $Y$ and the divisor $\Delta$ on the variety $R$. We get
$$
\mathop{\rm m}(Y\circ\Delta)\geqslant
\mathop{\rm m}(Y)+\mathop{\rm d}\nolimits_Q(B) \mathop{\rm mult}\nolimits_B Y^+\geqslant
\mathop{\rm m}(Y)+2\mathop{\rm mult}\nolimits_B Y^+.
$$
Recall that for $Z_R=Y$ the inequality (\ref{23.11.2017.1}) holds. It is easy to check that the minimum of the function $s+2t$ of the real variables $s,t$ on the set
$$
\{s+t\geqslant 2,\,\,s\geqslant 2t\}\subset{\mathbb R}^2_+
$$
is attained at the point $(\frac43,\frac23)$ and equal to $\frac83$. Therefore,
$$
\mathop{\rm m}(Y\circ\Delta)>\frac83 \mathop{\rm d}(Y\circ \Delta)
$$
($\Delta$ is a hyperplane section of $R$, so that
$\mathop{\rm d}(Y\circ\Delta)=\mathop{\rm d}(Y)$). The last inequality can be re-written in the form of the estimate
$$
\frac{\mathop{\rm mult}\nolimits_o}{\mathop{\rm deg}}(Y\circ\Delta)>\frac{16}{3M},
$$
which contradicts the inequality (\ref{04.11.2017.1}) for
$c=a=2$. We obtained a contradiction which proves the lemma. Q.E.D.\vspace{0.1cm}

By Lemmas 5.2 and 5.3, is the case 3.2 takes place, then on
$\Delta$ there is a prime divisor $Y$, satisfying the inequality
$$
\mathop{\rm m}(Y)+\mathop{\rm mult}\nolimits_B Y^+>2\mathop{\rm d}(Y),
$$
and moreover, $\mathop{\rm d}_Q (B)=2$ and
$\mathop{\rm m}(Y)>2\mathop{\rm mult}\nolimits_B Y^+$, so that
$\mathop{\rm m}(Y)>\frac43 \mathop{\rm d}(Y)$. In order to exclude that last option, we use the condition (R2.3). Denote the subvariety that is cut out on the quadric $\Delta_Q$ by the equation
$q_3|_{\Delta_Q}=0$, by the symbol $G$. It is easy to see that $G$ belongs to the family of varieties, which are irreducible, reduced and factorial by the condition (R2.3). Therefore, $G$ is a factorial complete intersection of type $2\cdot 3$ in ${\mathbb P}^{M-2}$. For that reason the kernel of the surjective restriction map
$$
\mathop{\rm res}\colon H^0({\mathbb P}^{M-2},
{\cal O}_{{\mathbb P}^{M-2}}(2))\to H^0(G,{\cal O}_G(2))
$$
is one-dimensional and is generated by the quadratic form $q_2|_{{\mathbb P}(\Pi)}$. The irreducible subvariety $B$ is cut out on the quadric
$\Delta_Q$ by a quadratic equation $\beta=0$, where
$\beta\not\in\langle q_2|_{{\mathbb P}(\Pi)}\rangle$. Therefore, the equation $\beta||_G=0$ defines an effective divisor on $G$.\vspace{0.1cm}

{\bf Lemma 5.4.} {\it The divisor $\{\beta|_G=0\}$ is irreducible and reduced.}\vspace{0.1cm}

{\bf Proof.} By the factoriality of the complete intersection
$G$, reducibility or non-reducedness of this divisor means that it is a sum of two hyperplane sections. Therefore, if the divisor $\{\beta|_G=0\}$ were reducible or non-reduced, for some linear forms
$l_1, l_2$ on ${\mathbb P}^{M-2}{\mathbb P}(\Pi)$ we would have had the equality
$$
\beta=l_1l_2+\lambda q_2|_{\Pi},
$$
where $\lambda\in{\mathbb C}$ is some constant. But then the divisor
$B\subset\Delta_Q$, given by the equation $\beta=0$, would have been reducible or non-reduced. Q.E.D. for the lemma.\vspace{0.1cm}

Therefore, $B\cap G=\{\beta|_G=0\}$ is an irreducible reduced subvariety of codimension 4 on $Q$, and moreover
$$
B\cap G\sim 6H^4_Q,
$$
that is, $\mathop{\rm d}_Q(B\cap G)=6$.\vspace{0.1cm}

Now let us come back to the prime divisor $Y$ on $\Delta$. If $Y^+$ does not contain $B$, then the inequality
$$
\mathop{\rm m}(Y)>2\mathop{\rm d}(Y)
$$
holds, which is excluded by the proof of Lemma 5.2 (where we proved that
$\mathop{\rm m}(Y)\leqslant \frac32 \mathop{\rm d}(Y)$). So
$Y+$ must contain $B$.\vspace{0.1cm}

Since $G=(T^+_2\circ \Delta_Q)=T^+_2\cap \Delta_Q$ does not contain the subvariety $B$, we conclude that $Y\not\subset T_2$. Therefore, the effective cycle $Y_*=(Y\circ T_2)$ of codimension 2 on $\Delta$ is well defined. The divisor $Y^+$ on $\Delta^+$ can, generally speaking, contain the subvariety $G\subset \Delta_Q$. For that reason, write down
$$
(Y^+\circ T^+_2)=bG+Y^+_*
$$
for some $b\in{\mathbb Z}_+$. Such a writing is possible, because
$(\Delta\circ T^+_2)$ by construction is precisely $G$. Since
$\mathop{\rm mult}\nolimits_{B\cap G} G=1$ and
$\mathop{\rm d}_Q(G)=3$, we get
$$
\mathop{\rm m}(Y_*)=3\mathop{\rm m}(Y)+3b,
$$
where we have
$$
\mathop{\rm mult}\nolimits_{B\cap G} Y^+_*\geqslant \mathop{\rm mult}\nolimits_B Y^+-b.
$$
Therefore, the inequality
$$
\mathop{\rm m}(Y_*)+3\mathop{\rm mult}\nolimits_{B\cap G} Y^+_*>6\mathop{\rm d}(Y)
$$
holds. Since $\mathop{\rm d}(Y_*)=2\mathop{\rm d}(Y)$, the last inequality can be re-written in the form
\begin{equation}\label{25.11.2017.1}
\mathop{\rm m}(Y_*)+3\mathop{\rm mult}\nolimits_{B\cap G} Y^+_*>6\mathop{\rm d}(Y_*).
\end{equation}
Besides, as we pointed out above,
$\mathop{\rm m}(Y)>\frac43\mathop{\rm d}(Y)$, which implies the inequality
\begin{equation}\label{25.11.2017.2}
\mathop{\rm m}(Y_*)>\frac32\cdot\frac43\mathop{\rm d}(Y_*)=2\mathop{\rm d}(Y_*),
\end{equation}
which will be used later.\vspace{0.1cm}

Consider again the irreducible reduced subvariety $B\cap G$. It can be viewed as a complete intersection of type $2\cdot 2\cdot 3$ in
${\mathbb P}(\Pi)={\mathbb P}^{M-2}$. Therefore $B\cap G$ is cut out in the scheme-theoretic sense by cubic hypersurfaces. Let
$$
|3H_{\Pi}-B\cap G|
$$
be the system of cubic hypersurfaces, containing $B\cap G$. Then
$\mathop{\rm Bs} |3H_{\Pi} - B\cap G|=B\cap G$. Furthermore, let
$$
C(|3H_{\Pi}-B\cap G|)
$$
be the corresponding linear system of cones in the projective space
$\overline{\Pi}\cong {\mathbb P}^{M-1}$, that is, a cubic hypersurface $W\in C(|3H_{\Pi}-B\cap G|)$ if and only if $W$ is a (cubic) cone with the vertex at the point $o$ and
$$
W^+\cap E_{\mathbb P}\in |3H_{\Pi}-B\cap G|.
$$
Obviously, $\mathop{\rm Bs} C(|3H_{\Pi}-B\cap G|)=C(B\cap G)$ is a cone with the vertex at the point $o$, consisting of lines $L\ni o$ such that
$L^+\cap E_{\mathbb P}\in B\cap G$.\vspace{0.1cm}

By the regularity condition $C(B\cap G)\not\subset\Delta$ (through the point
$o$ there are finitely many lines on $V$). The cone $C(B\cap G)$ is irreducible, so that for the restriction
$$
C(|3H_{\Pi}-B\cap G|)_{\Delta}
$$
of the linear system $C(|3H_{\Pi}-B\cap G|)$ on $\Delta$ we have:
$$
\mathop{\rm codim}(\mathop{\rm Bs} C(|3H_{\Pi}-B\cap G|)_{\Delta}\subset\Delta)=3
$$
and for that reason for a general divisor $W\in C(|3H_{\Pi}-B\cap G|)$ none of the components of the support of the cycle $Y_*$ of codimension 2 in $\Delta$ is contained in $W$, so that the effective cycle $(Y_*\circ W)$ of their scheme-theoretic intersection is well defined. This is a cycle of codimension 3 on $\Delta$. For this cycle we get the estimate
$$
\mathop{\rm mult}\nolimits_o(Y_*\circ W)\geqslant 3\mathop{\rm mult}\nolimits_o Y_*+
12\mathop{\rm mult}\nolimits_{B\cap G} Y^+_*
$$
(since $\mathop{\rm deg} B\cap G=12$), which we re-write in the form
$$
\mathop{\rm mult}\nolimits_o(Y_*\circ W)\geqslant 6(\mathop{\rm m} (Y_*)+
2\mathop{\rm mult}\nolimits_{B\cap G} Y^+_*).
$$
Obviously, $\mathop{\rm d}(Y_*\circ W)=3\mathop{\rm d}(Y_*)$. Taking into account that the minimum of the function $s+2t$ of real variables $s,t$ on the set
$$
\{s+3t\geqslant 3,\,\, s\geqslant 2\}\subset {\mathbb R}^2_+
$$
is attained at the point $(2,\frac13)$ and is equal to $\frac83$, we get by the inequalities (\ref{25.11.2017.1}) and (\ref{25.11.2017.2}):
$$
\mathop{\rm mult}\nolimits_o(Y_*\circ W)>6\cdot \frac83 \mathop{\rm d}(Y_*)=\frac{16}{3}
\mathop{\rm d}(Y_*\circ W),
$$
that is, the inequality
$$
\frac{\mathop{\rm mult}\nolimits_o}{\mathop{\rm deg}} (Y_*\circ W)>\frac{16}{3M}
$$
holds. This contradicts the inequality (\ref{04.11.2017.1}) for
$c=2$, $a=3$ when $M\geqslant 16$. We obtained a contradiction which completes the exclusion of the case 3.2.\vspace{0.3cm}


{\bf 5.3. Exclusion of all remaining cases.} Assume that one of the following four cases takes place:
$$
3.1,\quad 3.0,\quad 2.1,\quad 2.0.
$$
Let us consider the self-intersection $Z$ of the mobile system $\Sigma$.\vspace{0.1cm}

{\bf Lemma 5.5.} {\it The following inequality holds:}
$$
\mathop{\rm m}(Z)\geqslant 2\mathop{\rm mult}\nolimits_B Z^+.
$$

{\bf Proof.} In the notations of \S 4 let $P\in {\cal L}$ be a general subspace of maximal dimension on the quadric $Q$, $Z_P=(Z^+\circ P)$ an effective cycle of codimension 2. Assume the converse:
$$
\mathop{\rm m}(Z)<2\mathop{\rm mult}\nolimits_B Z^+.
$$
Then $Z_P$ is an effective cycle of degree
$\mathop{\rm m}(Z)$ (recall that
$\mathop{\rm m}(Z)=\mathop{\rm d}_Q (Z^+\circ Q)$), which satisfies the inequality
$$
\mathop{\rm deg} Z_P<2\mathop{\rm mult}\nolimits_{B\cap P} Z_P.
$$
Let $p,q\in B\cap P$ be a pair of distinct points and
$[p,q]\subset P$ the line through them. Furthermore, let
$\Theta\supset [p,q]$ be a two-dimensional plane of general position in
$P$, containing the line $[p,q]$. By the symbol $|Z_P|$ we denote the support of the cycle $Z_P$. If the intersection $\Theta\cap |Z_P|$ were zero-dimensional, we would have got the following chain of equalities and inequalities, where all intersection numbers are meant to be in the projective space $P$:
$$
\mathop{\rm deg} Z_P=(Z_P\cdot \Theta)=\sum_{x\in \Theta\cap |Z_P|} (Z_P\cdot\Theta)_x\geqslant
$$
$$
\geqslant (Z_P\cdot\Theta)_p+(Z_P\cdot\Theta)_q\geqslant
\mathop{\rm mult}\nolimits_p Z_P+\mathop{\rm mult}\nolimits_q Z_P\geqslant 2\mathop{\rm mult}\nolimits_{B\cap P} Z_P>
\mathop{\rm deg} Z_P,
$$
which is impossible. Therefore, the set $\Theta\cap |Z_P|$ is positive-dimensional. Since the plane $\Theta$ is arbitrary, we conclude that $[p,q]\subset |Z_P|$. Therefore, the support $|Z_Q|$ of the cycle
$Z_Q=(Z^+\circ Q)$ contains the secant subvariety
$\mathop{\rm Sec} (B\subset Q)$. However,
$\mathop{\rm codim} (|Z_Q|\subset Q)=2$, whereas by assumption we are in one of the four cases when
$$
\mathop{\rm codim}(\mathop{\rm Sec} (B\subset Q)\subset Q)\leqslant 1.
$$
This contradiction proves the lemma. Q.E.D.\vspace{0.1cm}

Now let us use the technique of counting multiplicities (Subsection 3.3). In the notations of Subsection 3.3 Lemma 5.5 claims that $m_1\geqslant 2m_2$.\vspace{0.1cm}

{\bf Lemma 5.6.} {\it The following inequality holds:}
$$
m_1>8n^2.
$$

{\bf Proof.} Assume the converse:
$m_1\leqslant 8n^2$. Then $m_2\leqslant 4n^2$. From the inequality
(\ref{15.11.2017.4}) we get:
$$
(8p_1+4(\Sigma^+_0+\Sigma_1))(p_1+\Sigma_0^++\Sigma_1+\Sigma_2)>
(3p_1+3\Sigma_0^++2\Sigma_1+\Sigma_2)^2.
$$
It is easy to bring the last inequality to the following form:
$$
0>(p_1-\Sigma_2)^2+(6p_1+5\Sigma_0^++4\Sigma_1+2\Sigma_2)\Sigma_0^+.
$$
The contradiction proves the lemma. Q.E.D.\vspace{0.1cm}

Thus in each of the four cases under consideration the self-intersection
$Z$ satisfies the inequality
$$
\frac{\mathop{\rm mult}\nolimits_o}{\mathop{\rm deg}}Z=\frac{2}{M}\cdot
\frac{\mathop{\rm m}(Z)}{\mathop{\rm d}(Z)}>\frac{4}{M}.
$$
This contradicts Proposition 2.4 (see Remark 2.1 (i)). The cases 3.1, 3.0, 2.1 and 2.0 are excluded.\vspace{0.1cm}

The only remaining case is the case 2.2. Assume that this case takes place. Recall that $B$ in this case is a section of the quadric $Q$ by a linear subspace of codimension 2 in $E_{\mathbb P}\cong{\mathbb P}^M$. Let
$\Delta$ be the section of the hypersurface $V$ by the uniquely determined subspace of codimension 2 in ${\mathbb P}$, such that
$\Delta^+\cap Q=B$. Then
$\mathop{\rm m}(\Delta)=\mathop{\rm d}(\Delta)=\mathop{\rm mult}\nolimits_B\Delta^+=1$. Write down
$$
Z=a\Delta+Z_1,
$$
where $a\in{\mathbb Z}_+$ and the cycle $Z_1$ does not contain $\Delta$ as a component. Since $\mathop{\rm d}(Z)=4n^2$, from Proposition 3.6 (and the equalities for $\Delta$, written above) we obtain:
$$
\mathop{\rm m}(Z_1)+\mathop{\rm mult}\nolimits_B Z^+_1>2\mathop{\rm d}(Z_1).
$$
By the linearity of this inequality we may assume that $Z_1=Y$ is an irreducible subvariety of codimension 2, and moreover, $Y\neq \Delta$. Consider a general hyperplane section $R\supset\Delta$. Since
$Y\not\subset R$, the effective cycle $(Y\circ R)$ of codimension 2 on $R$ is well defined. This cycle satisfies the inequality
$$
\mathop{\rm m}(Y\circ R)=\mathop{\rm m}(Y)+\mathop{\rm mult}\nolimits_B Y^+>
2\mathop{\rm d}(Y)=2\mathop{\rm d}(Y\circ R),
$$
which is equivalent to the estimate
$$
\frac{\mathop{\rm mult}\nolimits_o}{\mathop{\rm deg}} (Y\circ R)=
\frac{2\mathop{\rm m}(Y\circ R)}{M\mathop{\rm d}(Y\circ R)}>
\frac{4}{M},
$$
contradicting the inequality (\ref{04.11.2017.1}) for
$c=1$, $a=2$. This contradiction excludes the case 2.2 and completes the proof of Theorem 5.


\begin{flushleft}
Department of Mathematical Sciences,\\
The University of Liverpool
\end{flushleft}

\noindent{\it pukh@liverpool.ac.uk}


\begin{thebibliography}{99}
\bibitem{Pir2016} Pirutka A., Varieties that are not stably rational, zero-cycles and unramified cohomology, ArXiv: 1603.09261, 24 p.
\end{thebibliography}
\end{document}